\documentclass[]{interact}

\usepackage{epstopdf}
\usepackage[caption=false]{subfig}
\usepackage{enumitem}
\usepackage{algorithm}
\usepackage{graphicx}
\usepackage{subfig}
\usepackage{xcolor} 
\usepackage{algpseudocode}
\usepackage{float}
\usepackage{hyperref}
\usepackage{color}
\usepackage{multirow}
\usepackage[numbers,sort&compress]{natbib}
\bibpunct[, ]{[}{]}{,}{n}{,}{,}
\renewcommand\bibfont{\fontsize{10}{12}\selectfont}

\theoremstyle{plain}
\newtheorem{theorem}{Theorem}[section]

\newtheorem{conject}[theorem]{Conjecture}

\newenvironment{conclusions}
  {\FloatBarrier\par\medskip\noindent\textbf{Conclusions.}\ }
  {\par\medskip}

\theoremstyle{definition}
\newtheorem{definition}[theorem]{Definition}

\theoremstyle{remark}
\newtheorem{remark}{Remark}

\usepackage{placeins}

\begin{document}


\title{Numerical Design of Optimized First-Order Algorithms}

\author{
\name{Yassine Kamri\textsuperscript{a}\thanks{CONTACT Yassine Kamri. Email: yassinekamri1996@gmail.com}, Julien M. Hendrickx\textsuperscript{a} and François Glineur\textsuperscript{a}}
\affil{\textsuperscript{a} UCLouvain, ICTM, INMA, Avenue Georges Lemaître, 4, 1348 Louvain-la-Neuve, Belgium}
}

\maketitle

\begin{abstract}
We derive several numerical methods for designing optimized first-order algorithms in unconstrained convex optimization settings. Our methods are based on the Performance Estimation Problem (PEP) framework, which casts the worst-case analysis of optimization algorithms as an optimization problem itself. We benchmark our methods against existing approaches in the literature on the task of optimizing the step sizes of memoryless gradient descent (which uses only the current gradient for updates) over the class of smooth convex functions. We then apply our methods to numerically tune the step sizes of several memoryless and full (i.e., using all past gradient information for updates) fixed-step first-order algorithms, namely coordinate descent, inexact gradient descent, and cyclic gradient descent, in the context of linear convergence. In all cases, we report accelerated convergence rates compared to those of classical algorithms.
\end{abstract}


\section{Introduction}
First-order algorithms have become some of the methods of choice for solving modern optimization problems, which are increasingly complex and large-scale, due to their simplicity and low computational cost. However, their effectiveness heavily depends on the choice of step sizes. As a result, finding step sizes that lead to more efficient first-order algorithms has become a major focus of research. The first results on accelerating the convergence of first-order methods date back to the work of Nemirovsky and Polyak~\cite{nemirovsky1984chebyshev,polyak1964methods}, where accelerated methods for quadratic problems were presented. Beyond the quadratic case, efforts to accelerate the convergence of first-order algorithms began with the work of Nemirovsky and Yudin on information-based complexity~\cite{nemirovsky1983complexity}. Nemirovsky also proposed the first accelerated first-order method for smooth convex optimization, which required a two-dimensional subspace search. Nesterov later resolved this limitation by introducing a computationally efficient accelerated gradient method for smooth convex functions~\cite{nesterov1983FGM}.
 Since then, acceleration schemes have been extended to a wide range of problems, including proximal gradient methods~\cite{beck2009fista,lin2015universal,monteiro2013accelerated}, random block coordinate descent~\cite{nesterov2012Coords,fercoq2015coord}, and stochastic gradient descent~\cite{shalev2014accelerated,allen2017katyusha}. The recent development of the PEP framework~\cite{drori2014perf,taylor2017smooth}, which provides a way to automatically compute upper bounds on the convergence of fixed-step first-order methods in various settings, has simplified the design of optimized first-order methods. This has led to the discovery of new accelerated gradient methods, some of which are information-theoretically optimal for their respective problem classes~\cite{kim2016OGM,taylor2023ITEM}. Despite these successes, using PEPs to design first-order methods results in nonconvex problem formulations. Consequently, existing approaches have largely relied on heuristic step-size schedules~\cite{GrimmerLongsteps,grimmer2024composing,Altschuler2023I,Altschuler2024II}, problem-specific relaxations~\cite{barre2023principled,taylor2023ITEM,drori2020GFM,kim2016OGM,goujaud2022optimal,park2022exact} where analytical solutions were found, or computationally intensive branch-and-bound methods (BNB-PEP)~\cite{dasgupta2024BNB,jang2023optista}. In this paper, our objective is to develop methods for the automatic, computer-assisted numerical design of optimized fixed-step first-order optimization algorithms. We apply these methods to several convex optimization settings and show that our optimized algorithms can achieve faster convergence rates than classical first-order methods in these settings. This paper is organized as follows:

\medskip

\noindent \textbf{Organization and contributions.} In Section~\ref{sec:design_prob}, we start by defining a generic approach to the design of fixed-step first-order methods. We then present several numerical approaches for automatically designing optimized first-order methods in Section~\ref{sec:design_methods}. In Section~\ref{sec:bench}, we compare the efficiency of our methods with the BNB-PEP approach on the benchmark task of designing an optimized fixed-step first-order memoryless gradient descent algorithm over the class of smooth convex functions. Finally, in Section~\ref{sec:opt_algo}, we apply our methodologies to numerically design optimized fixed-step first-order algorithms in several different settings, namely for block coordinate-wise algorithms, inexact gradient descent and a cyclic gradient descent algorithm in the setting of linear convergence of the gradient. We provide improved numerical upper bounds on convergence rates for these algorithms compared to existing fixed-step first-order algorithms. Furthermore, our numerical results suggest that recent findings on the accelerated convergence of first-order methods over smooth (strongly) convex functions~\cite{Altschuler2023I,Altschuler2024II,grimmer2024composing} without using Nesterov momentum-based acceleration schemes can be extended to several other settings.

\section{Problem Formulation}\label{sec:design_prob}

In this section, we define the problem of designing optimized fixed-step first-order algorithms, for certain classes of functions. To do so, we first show how to cast the worst-case performance analysis of a fixed-step first-order algorithm as an optimization problem itself which was first done in~\cite{drori2014perf} for gradient descent over smooth (strongly) convex functions. fixed-step first-order algorithms are defined as follows:
\begin{definition}{(Fixed-step first-order method).}
An algorithm \( \mathcal{M} \) is called a \emph{fixed-step first-order method} if its iterates are computed as
\begin{equation}\label{eq:fixed_step_alg}
    x_i = x_0 - \sum_{k=0}^{i-1} \alpha_{i,k} g_k,
\end{equation}
with fixed scalar coefficients \( \alpha_{i,k} \). We denote this method by \( \mathcal{M}_{\alpha} \).
\end{definition}
\noindent Consider a class of convex functions \( \mathcal{F}(\mathbb{R}^d) \) and a first-order fixed-step algorithm  \( \mathcal{M}_\alpha \), the worst-case performance of algorithm \( \mathcal{M}_\alpha \) over the functional class \( \mathcal{F}(\mathbb{R}^d) \) is given by the optimal solution of the following infinite-dimensional optimization problem:

\begin{equation}\label{i-PEP}
\begin{aligned}
    \sup_{f,\, x_0,\, \ldots,\, x_N,\, x_*} \quad & \mathcal{P}(f, x_0, \ldots, x_N, x_*) \\
    \text{subject to} \quad 
        & f \in \mathcal{F}(\mathbb{R}^d),\\
        & x_* \text{ is optimal for } f,\\
        & x_1, \ldots, x_N \text{ are generated from } x_0 \text{ by method } \mathcal{M}_\alpha,\\
        & \mathcal{I}(f, x_0, \ldots, x_N, x_*) \leq R,
\end{aligned}
\end{equation}
where $\mathcal{P}$ is the performance criterion for instance the final objective accuracy $\mathcal{P}(f, x_0, \ldots, x_N, x_*) = f(x_N) - f(x_*)$, and  $\mathcal{I}(f, x_0, \ldots, x_N, x_*) \leq R$ an initial condition to ensure that the worst-case of the algorithm is bounded, for instance the initial distance to a solution $\mathcal{I}(f, x_0, \ldots, x_N, x_*) = \|x_0 - x_*\|$. Infinite-dimensional problems such as~\eqref{i-PEP} are typically challenging to solve. To address this, Taylor et al.~\cite{taylor2017smooth} introduced the concept of \( \mathcal{F}(\mathbb{R}^d) \)-interpolability of finite sets:
\begin{definition}[\( \mathcal{F}(\mathbb{R}^d) \)-interpolability]\label{def:F_interp}
Let \( I \) be an index set, and consider the set of triples \( \mathcal{S} = \{(x_i, g_i, f_i)\}_{i \in I} \), where \( x_i, g_i \in \mathbb{R}^d \) and \( f_i \in \mathbb{R} \) for all \( i \in I \). The set \( \mathcal{S} \) is \( \mathcal{F}(\mathbb{R}^d) \)-interpolable if and only if there exists a function \( f \in \mathcal{F}(\mathbb{R}^d) \) such that $g_i \in \partial f(x_i)$ and $f(x_i) = f_i, \; \forall i \in I$.
\end{definition}
\noindent which allows Problem~\eqref{i-PEP} to be reformulated as the following finite-dimensional Problem~\eqref{i-PEP}:

\begin{equation}\label{PEP-generic-finite}
\begin{aligned} 
& \underset{\mathcal{S}_N = \{(x_i,g_i,f_i)\}_{i \in I}}{\text{sup}}
& & \mathcal{P}(\mathcal{S}_N) \\ 
& \text{subject to}
& & \mathcal{I}(\mathcal{S}_N) \leq R ,\\
& & & \mathcal{S}_{N} \text{\;is\;} \mathcal{F}(\mathbb{R}^d)\text{-interpolable} ,\\
& & &  x_0,\dots,x_N \; \text{generated by} \; \mathcal{M}_{\alpha} ,\\
& & & g_* = 0.
\end{aligned}
\end{equation}
Using this, we can conceptually define the problem of finding an optimal fixed-step first-order algorithm over the class $\mathcal{F}(\mathbb{R}^d)$ (i.e, the algorithm with the lowest possible worst-case over the class $\mathcal{F}(\mathbb{R}^d)$) as solving the following min-max problem.
\begin{equation}\label{PEP-generic-opt}
\begin{aligned} 
& \underset{\alpha}{\text{inf}} \;\; \underset{\mathcal{S}_N = \{(x_i,g_i,f_i)\}_{i \in I}}{\text{sup}}
& & \mathcal{P}(\mathcal{S}_N) \\ 
& \text{subject to}
& & \mathcal{I}(\mathcal{S}_N) \leq R,\\
& & & \mathcal{S}_N \text{ is } \mathcal{F}(\mathbb{R}^d)\text{-interpolable},\\
& & & x_0,\dots,x_N \text{ generated by } \mathcal{M}_{\alpha},\\
& & & g^* = 0.
\end{aligned}
\end{equation}
For certain types of functional classes, the PEP defined by Problem~\eqref{PEP-generic-finite} can be reformulated as a convex conic optimization problem, thereby simplifying the algorithm design problem~\eqref{PEP-generic-opt}. We call these functional classes SDP-representable classes of functions: 

\begin{definition}[SDP-representable classes of functions]\label{def:class_func_sdp} 
Given a class of functions \( \mathcal{F}(\mathbb{R}^d) \), we say that \( \mathcal{F}(\mathbb{R}^d) \) is SDP-representable if and only if, for any set \( \mathcal{S}_{N} = \{(x_i,g_i,f_i)\}_{i \in I} \) with \( I = \{0,\dots,N,*\} \), where \( \{x_i\}_{i \in \{0,\dots,N\}} \) are the iterates generated by the fixed-step first-order algorithm \( \mathcal{M}_{\alpha}\) defined by~\eqref{eq:fixed_step_alg}, the functional class \( \mathcal{F}(\mathbb{R}^d) \) admits necessary interpolation conditions such that if the set \( \mathcal{S}_{N} \) is \( \mathcal{F}(\mathbb{R}^d) \)-interpolable, then:  
\begin{equation}\label{eq:gen_interp_conds}
    F^\top (\mu_i u_i - \mu_j u_j) \geq \operatorname{Tr}(A_{i,j}(\alpha) G), \; \forall i,j \in I,
\end{equation}
where $G = P^\top P \in \mathbb{S}^{N+2}$, $P = [g_0,g_1,\dots,g_N,x_0]$, \( F^\top = [f_0,f_1,\dots,f_N,f_*] \), $u_i = e_{i+1} \in \mathbb{R}^{N+2}, \text{ where } \{e_i\}_{i \in \{1,\ldots,N+2\}} \text{ is the canonical basis of } \mathbb{R}^{N+2}$, \( \mu_i \) for all \( i \in I \) are scalars (usually $\mu_i$ = $\mu_j = 1$) and \( A_{i,j}(\alpha) \in \mathbb{R}^{(N+2)\times(N+2)}\) .  
\end{definition}

\begin{theorem}\label{th:generic_PEP}
Consider an SDP-representable functional class \( \mathcal{F}(\mathbb{R}^d) \), the fixed-step first-order method \( \mathcal{M}_{\alpha}\), a performance criterion of the form \( \mathcal{P}(F, G) = b^{\top} F + \operatorname{Tr}(C G) \), and an initial condition of the form \( \mathcal{I}(G,F)^2 = \operatorname{Tr}(A_R G) + \|v^\top F\|^2 \). For any \( d \in \mathbb{N} \), a valid upper bound on the performance after \( N \) steps of the fixed-step first-order algorithm \( \mathcal{M}_{\alpha}\) over the class of functions \( \mathcal{F}(\mathbb{R}^d) \) is given by:  
\begin{equation}\tag{SDP-PEP-Generic}\label{eq: SDP-PEP-Generic}
\begin{aligned}
      w^{\text{SDP}}(R, \alpha, N, b, C,A_R ) = &\sup_{G \in \mathbb{S}^{N+2}, \; f \in \mathbb{R}^{N+1}} \; b^{\top} F + \operatorname{Tr}(C G) \\
    &\text{such that:} \\
    & F^\top (\mu_i u_i - \mu_j u_j) \geq \operatorname{Tr}(A_{i,j}(\alpha) G), \; \forall i,j \in I,\\
    &  \operatorname{Tr}(A_R G) + \|v^\top F\|^2 - R^2 \leq 0,\\
    & G \succeq 0,
\end{aligned}
\end{equation}
where $G = P^\top P \in \mathbb{S}^{N+2}$, $P = [g_0,g_1,\dots,g_N,x_0]$, \( F^\top = [f_0,f_1,\dots,f_N,f_*] \), $u_i = e_{i+1} \in \mathbb{R}^{N+2}, \text{ where } \{e_i\}_{i \in \{1,\ldots,N+2\}} \text{ is the canonical basis of } \mathbb{R}^{N+2}$, \( \mu_i \) for all \( i \in I \) are scalars (usually $\mu_i$ = $\mu_j = 1$) and \( A_{i,j}(\alpha) \in \mathbb{R}^{(N+2)\times(N+2)}\).   
\end{theorem}
\begin{proof}
  Taylor et al.~\cite{taylor2017smooth} provided necessary and sufficient interpolation conditions for the classes of $L$-smooth $\mu$-convex functions, showing that these classes are SDP-representable, and constructed a convex SDP formulation for the corresponding Performance Estimation Problems. For a generic construction of the convex Performance Estimation Problem for SDP-representable classes of functions, we refer the reader to~\cite[Chapter~5, Section~5.2]{kamri2025PhD}.
\end{proof}

\begin{remark}
Note that, since we only require the interpolation conditions to be necessary (but not sufficient), the SDP formulation presented in the previous theorem may be a relaxation of the original Problem~\ref{PEP-generic-finite}. As a result, it provides a valid upper bound on the worst-case performance of the algorithm over the functional class, but not necessarily the exact worst-case. However, if the interpolation conditions are shown to be both necessary and sufficient, then the SDP is an exact reformulation of Problem~\ref{PEP-generic-finite}, and its optimal value corresponds to the exact worst-case performance of the algorithm over the class.
\end{remark}
\noindent Using the generic convex PEP framework given in Theorem~\ref{th:generic_PEP}, we can now formulate the problem of designing optimized fixed-step first-order algorithms for any SDP-representable functional class as the following min-max problem:  

\begin{equation*}
\begin{aligned} 
& \underset{\alpha_{i,k}}{\text{inf}} \; \underset{F,G}{\text{sup}}
& & b^{\top}F + \mathrm{Tr}(CG)\\ 
& \text{subject to}
& &  \mathrm{Tr}(A_RG)\leq R^2,\\
& & & F^\top (\mu_i u_i - \mu_j u_j) \geq \mathrm{Tr}(A_{i,j}(\alpha)G), \; \forall i,j \in I,\\
& & & G \succeq 0.
\end{aligned}
\end{equation*}
\begin{remark}
Note that since the convex PEP formulation in Theorem~\ref{th:generic_PEP} may be a relaxation, this, in turn, can induce a relaxation in the algorithm design problem formulation derived above. As a result, solving it does not guarantee the optimality of the algorithm for the considered functional class.
\end{remark}
\noindent Since min-max optimization problems are typically hard to solve, we can further simplify the formulation by transforming it into a single minimization problem. This can be achieved by computing the Lagrangian dual of Problem~\eqref{eq: SDP-PEP-Generic} (i.e., dualizing the inner problem). The dual of Problem~\eqref{eq: SDP-PEP-Generic} takes the form:  
\begin{equation}\tag{dual-PEP-generic}\label{dual-PEP-generic}
\begin{aligned} 
w^{\text{dual}}(\alpha) =\; & \underset{\tau,\lambda_{i,j}}{\text{inf}}
& & \tau R^2\\ 
& \text{subject to}
& &  \tau A_R - C + \sum_{i,j \in I} \lambda_{i,j} A_{i,j}(\alpha) \succeq 0 ,\\
& & & b - \sum_{i,j} \lambda_{i,j} (\mu_j u_j - \mu_i u_i) = 0,\\
& & &  \lambda_{i,j} \geq 0, \; \forall i,j \in I,\\
& & &  \tau \geq 0 .\\
\end{aligned}
\end{equation}
\begin{remark}
Note that if strong duality does not hold, then using the Lagrangian dual of Problem~\eqref{eq: SDP-PEP-Generic} still yields a valid upper bound on the worst-case performance of the algorithm. However, this may further relax the corresponding algorithm design problem. Experimentally, strong duality holds for every explicit PEP formulation we consider in this paper.
\end{remark}
\noindent Thus, the problem of designing fixed-step first-order algorithms over an SDP-representable class of functions \( \mathcal{F}(\mathbb{R}^d) \) takes the simpler form:  
\begin{equation}\tag{Algo-design-generic}\label{algo-design-generic}
\begin{aligned} 
& \underset{\alpha_{i,k},\tau,\lambda_{i,j}}{\text{inf}}
& & \tau R^2\\ 
& \text{subject to}
& &  \mathcal{A}(\tau,\lambda,\alpha) \succeq 0 ,\\
& & & b - \sum_{i,j} \lambda_{i,j} (\mu_j u_j - \mu_i u_i) = 0,\\
& & &  \lambda_{i,j} \geq 0, \; \forall i,j \in I ,\\
& & &  \tau \geq 0 ,\\
\end{aligned}
\end{equation}
where \( \mathcal{A}(\tau,\lambda,\alpha)  = \tau A_R - C + \sum_{i,j \in I} \lambda_{i,j} A_{i,j}(\alpha) \). We obtain a nonconvex conic formulation of the algorithm design problem for any SDP-representable class of functions, similar to the formulations derived in~\cite{kim2016OGM,taylor2017smooth,dasgupta2024BNB} for the class of convex smooth functions.  
\begin{remark}
The matrix \( \mathcal{A}(\tau,\lambda,\alpha) \) is linear in the variables \( \lambda \) when the \( \alpha \) are fixed. However, it is nonlinear in both \( \alpha \) and \( \lambda \) together and potentially nonlinear in the variables \( \alpha \) when the \( \lambda \) are fixed. This implies that the constraint \( \mathcal{A}(\tau,\lambda,\alpha) \succeq 0 \) is a nonlinear matrix inequality. In general, optimization problems involving nonlinear matrix inequalities are NP-hard to solve \cite{toker1995np}. Matrix \( \mathcal{A}(\tau,\lambda,\alpha) \) is always linear in $\tau$.
\end{remark}

\section{Numerical Design Methods for First-Order Optimization Algorithms} \label{sec:design_methods}
We present several methods for solving the problem of designing optimized first-order methods, namely: an alternating minimization approach, a first-order method, and a sequential linearization method.

\noindent \textbf{Alternating Minimization (AM).} Although Problem~\eqref{algo-design-generic} is a nonlinear SDP, in many cases of interest, it exhibits a bilinear structure. Specifically, in these cases, the matrix \( \mathcal{A}(\tau,\lambda,\alpha) \) is also linear with respect to the variables \( \alpha \) when the \( \lambda \) are fixed. This is notably the case when designing fixed-step first-order algorithms for the class of \( L \)-smooth convex functions \( \mathcal{F}_{0,L}(\mathbb{R}^d) \) \cite{kim2016OGM,kim2021OGM-G}. In the case of designing algorithms over the class of \( L \)-smooth \( \mu \)-strongly convex functions \( \mathcal{F}_{\mu,L}(\mathbb{R}^d) \) with \( 0 \leq \mu < L < +\infty \), the problem can be reduced to designing algorithm over the class of $L-\mu$-smooth convex functions \( \mathcal{F}_{0,L-\mu}(\mathbb{R}^d) \)  as shown in \cite{taylor2023ITEM}. In these instances, when either the variables \( \alpha \) or \( \lambda \) are fixed, Problem~\eqref{algo-design-generic} reduces to a linear SDP. Therefore, it is possible to alternately optimize with respect to \( \lambda \) and \( \alpha \) by solving linear SDPs, following the approach in~\cite{Boyd1994BMI}. Let us now define the subproblem when the \( \lambda \) are fixed:  

\begin{equation}\tag{Opt-steps}\label{Opt-steps}
\begin{aligned} 
& \underset{\alpha_{i,k},\tau}{\text{inf}}
& & \tau R^2\\ 
& \text{subject to}
& & \mathcal{A}(\tau,\lambda,\alpha) \succeq 0 ,\\
& & & b - \sum_{i,j} \lambda_{i,j} (\mu_j u_j - \mu_i u_i) = 0,\\
& & &  \tau \geq 0 .\\
\end{aligned}
\end{equation}
Note that variable $\tau$ is optimized at both steps of our algorithm. We can now define our alternating minimization algorithm for the design of first-order methods:
\begin{algorithm}[H]
\caption{Alternating minimization for optimizing first-order algorithms (AM)}
\begin{algorithmic}[1]
    \State \textbf{Input:} Initial step-sizes $\alpha^{(0)}_{i,k}$, Number of iterations $T$
    \State \textbf{Output:} $\alpha^{(T)}_{i,k}$
    \For{$t = 1$ to $T$}
        \State Compute $\lambda_{i,j}^{(t)}$ by solving Problem~\eqref{dual-PEP-generic} for $\mathcal{M}_{\alpha^{(t-1)}}$
        \State Compute $\alpha_{i,k}^{(t)}$ by solving problem \eqref{Opt-steps} with $\lambda_{i,j} = \lambda_{i,j}^{(t)}$
    \EndFor
\end{algorithmic}
\end{algorithm}

\begin{remark}
Note that when strong duality holds, the variables \( \lambda_{i,j} \) in Problem~\eqref{dual-PEP-generic} characterize the proof of the upper bound on the performance of the algorithm \( \mathcal{M}_{\alpha} \). Informally, we obtain the upper bound by multiplying each constraint in Problem~\eqref{eq: SDP-PEP-Generic} by its corresponding dual variable \( \lambda_{i,j} \) and summing all the resulting inequalities. For a more detailed explanation of PEP-assisted convergence proofs, we refer the reader to~\cite{Goujaud2023proofs}. Thus, our alternating minimization approach can be understood as alternatively finding the optimal upper bound on the performance of a given algorithm along with a valid proof with respect to the set of constraints used in Problem~\eqref{eq: SDP-PEP-Generic}, and then finding the best algorithm such that the proof of the previous upper bound remains valid for the new algorithm.  
\end{remark}

\begin{remark}
\noindent Note that although our alternating minimization approach does yield better step-sizes compared to the initial step-sizes, we have no guarantee that it converges to a global minimum nor even to a local one. Indeed consider the simple example of optimizing two steps of gradient descent over $L$-smooth convex functions. In the next figure, we provide a contour plot for step sizes between $0$ and $\frac{2}{L}$ of the worst-case performance obtained by PEP and highlight the critical points of the worst-case over the space of step sizes (one global minimum and one local minimum), we also place on this figure the successive iterates of alternating minimization initialized with the classical step sizes $\frac{1}{L}$. Figure~\ref{fig:contour_plot} shows that alternating minimization gets stuck at a point of non-smoothness at the boundary between two regions and does not converge to either minimum.

\begin{figure}[!htbp]
    \centering
    \includegraphics[width=0.7\linewidth]{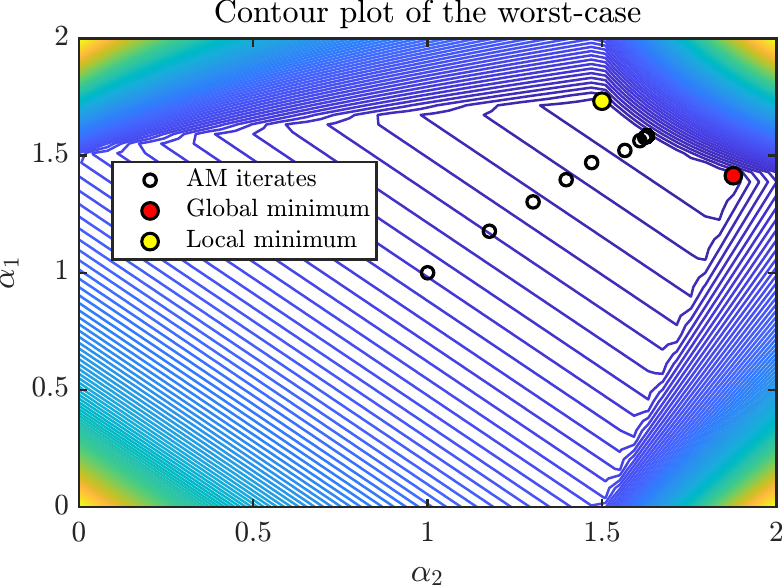}  
    \caption{AM iterates are not guaranteed to converge to a minimum } 
    \label{fig:contour_plot}
\end{figure}
\end{remark}

\noindent \textbf{First-order method (FOM).} Our second approach is to apply a first-order method to the function \( w^{\text{dual}}(\alpha) \) directly. This requires computing the gradients, when well-defined, of \( w^{\text{dual}} \) with respect to the step sizes \( \alpha_{i,k} \). Since the step sizes \( \alpha_{i,k} \) influence \( w^{\text{dual}}(\alpha) \) only through the vectors \( h_i \) in the matrices \( A_{i,j}(\alpha) \) (see Problem~\ref{dual-PEP-generic}), computing the gradients of \( w^{\text{dual}} \) with respect to the step sizes requires characterizing the variations of \( w^{\text{dual}} \) under small perturbations of the matrix  
\[
\tau A_R - C + \sum_{i,j \in I} \lambda_{i,j} A_{i,j}(\alpha).
\]
To this end, we consider for any symmetric matrix $B$, the perturbed problem:  

\begin{equation}\tag{Perturbed-dual-PEP}\label{Perturbed-dual-PEP}
\begin{aligned} 
\phi(B,\alpha) =\; & \underset{\tau,\lambda_{i,j}}{\text{inf}}
& & \tau R^2 \\ 
& \text{subject to}
& &  \tau A_R - C + \sum_{i,j \in I} \lambda_{i,j} A_{i,j}(\alpha) + B \succeq 0 ,\\
& & & b - \sum_{i,j} \lambda_{i,j} (\mu_j u_j - \mu_i u_i) = 0,\\
& & &  \lambda_{i,j} \geq 0 ,\\
& & &  \tau \geq 0 .\\
\end{aligned}
\end{equation}
\begin{theorem}{\cite[Proposition 5.80]{BonnansShapiro2000}}
The function $\phi(.,\alpha)$ is convex over the set of positive semidefinite matrices and:
\begin{equation*}
        \partial_B \phi(0,\{\alpha\}) = \{ -G:\; G\; \text{is a solution of Problem~\eqref{eq: SDP-PEP-Generic}}\}
\end{equation*}
\end{theorem}
\noindent This Theorem relies on the fact that the dual of Problem~\eqref{dual-PEP-generic} is Problem~\eqref{eq: SDP-PEP-Generic}.
\begin{remark}\label{diff_mat}
Note that, as a consequence of the previous theorem, the function \( \phi(.,\alpha) \) is differentiable at zero if and only if Problem~\eqref{eq: SDP-PEP-Generic} admits a unique solution.  
\end{remark}
\noindent Applying the chain rule, we have that:
\begin{theorem}\label{th:subgrad} if $w^{\text{dual}}$ is differentiable then the gradient is equal to the vector:
\begin{equation}\label{grad_PEP}
    \nabla w^{\text{dual}}(\alpha) = \left(\frac{\partial w^{\text{dual}}}{\partial \alpha_{1,0}}(\alpha),\dots, \frac{\partial w^{\text{dual}}}{\partial \alpha_{N,N-1}}(\alpha)    \right)
\end{equation}
\begin{equation*}
    \frac{\partial w^{\text{dual}}}{\partial \alpha_{m,n}}(\alpha) =  \operatorname{Tr} \left(-G \frac{\partial \mathcal{A}}{\partial \alpha_{m,n}}\right),
\end{equation*}
 where $G$ is the unique (see remark ~\ref{diff_mat}) solution of Problem~\eqref{eq: SDP-PEP-Generic}, $\mathcal{A}$ is defined by:
 \begin{equation*}
     \mathcal{A} = \tau A_R - C + \sum_{i,j \in I} \lambda_{i,j} A_{i,j}(\alpha)
 \end{equation*}
 and $\frac{\partial \mathcal{A}}{\partial \alpha_{m,n}}$ is the matrix such that:
 \begin{equation*}
    \left( \frac{\partial \mathcal{A}}{\partial \alpha_{m,n}}\right)_{i,j} =  \frac{\partial \mathcal{A}_{i,j}}{\partial \alpha_{m,n}}
 \end{equation*}
\end{theorem}
\begin{remark}
In practice, we use the vector computed in equation~\eqref{grad_PEP}, even though we cannot guarantee that it is a gradient everywhere or even an ascent direction. Numerical observations show that usually all possible step sizes choices for an algorithm can be divided into finitely many regions within which the PEP admits a unique solution, making $w^{\text{dual}}$ differentiable (see remark~\ref{diff_mat}) on these regions. Multiple solutions for the PEP may occur only on the boundaries between these regions on which the function $w^{\text{dual}}$ can become non-differentiable (see the master theses~\cite{daccache2019performance,eloi2022worst}). Furthermore, our first-order method yields satisfactory results on the benchmark application presented in Section~\ref{sec:bench}.  
\end{remark}

\noindent We can now present our first-order approach for the design of optimized first-order algorithms:

\begin{algorithm}[H]
\caption{First-order method (FOM) for the design of first-order algorithms}
\begin{algorithmic}[1]
    \State \textbf{Input:} Initial step-sizes $\alpha^{(0)}_{i,k}$, Number of iterations $T$, a list of decreasing steps $\{\beta_t\}_{[1,\dots,T]}$
    \State \textbf{Output:} $\alpha^{(T)}_{i,k}$
    \For{$t = 1$ to $T$}
        \State Compute $\frac{\partial w^{\text{dual}}}{\partial \alpha_{m,n}}(\alpha^{(t-1)})$ using Theorem~\ref{th:subgrad}
        \State Update the step-sizes: $\alpha^{(t)}_{m,n} = \alpha^{(t-1)}_{m,n} - \beta_t \frac{\partial w^{\text{dual}}}{\partial \alpha_{m,n}}(\{\alpha^{(t-1)}_{i,k}\})$
    \EndFor
\end{algorithmic}
\end{algorithm}

\noindent \textbf{Sequential Linearization Method (SLM).} We now present a third method to optimize the step sizes of first-order algorithms, closely linked to sequential semidefinite programming~\cite{Correa2004SSP,Fares2002SSP}. The central idea is to iteratively linearize Problem~\eqref{algo-design-generic} around a current set of decision variables \( \{\tau, \lambda, \alpha\} \). To simplify notations, we define the vector of decision variables as \( x = \{\tau, \lambda_{i,j}, \alpha_{i,k}\} \) and denote by $x_i$ the $i$-th component of vector $x$. Using this notation, Problem~\eqref{algo-design-generic} takes the form:
\begin{equation}\label{opt_design_x}
\begin{aligned} 
& \underset{x}{\text{inf}}
& & c^{\top} x \\ 
& \text{subject to}
& &  \mathcal{A}(x) \succeq 0 ,\\
& & & b - d^{\top}x = 0,\\
& & &  a_i^{\top} x \geq 0, \; \forall i \in \{1,\dots,N^2+1\} .\\
\end{aligned}
\end{equation}
Note that all the constraints of the problem are linear except for $\mathcal{A}(x)$. Our method starts by choosing a first feasible point $x^{(0)}$ and iteratively generates a sequence of iterates $x^{(t)}$, by solving the tangent problem:
\begin{equation}\tag{OptAlgo-linear}\label{optPEP_lin}
\begin{aligned} 
& \underset{\Delta x}{\text{min}}
& & c^{\top} \Delta x \\ 
& \text{subject to}
& &  \mathcal{A}(x^{(t)}) + \sum \frac{\partial \mathcal{A}}{\partial x_i}(x^{(t)}) \Delta x_i \succeq 0 ,\\
& & & b - d^{\top}(x^{(t)} + \Delta x)= 0,\\
& & &  a_i^{\top} (x^{(t)} + \Delta x) \geq 0 ,\\
\end{aligned}
\end{equation}
where $[\frac{\partial \mathcal{A}}{\partial x_i}(x^{(t)})]_{m,n} = \frac{\partial [\mathcal{A}(x^{(t)})]_{m,n} }{\partial x_i}$.
\newline
\begin{remark}
Our approach is very similar to sequential semidefinite programming (SQP); the only difference is that, typically, the tangent problem considered in SQP methods includes a quadratic regularization term \( \frac{1}{2} \Delta x^{\top} \tilde{H}_t \Delta x \) in the objective, where the matrix \( \tilde{H}_t \) is a positive semidefinite approximation of the Hessian matrix of the Lagrangian of the problem. In our case, the Hessian of the Lagrangian of Problem~\eqref{opt_design_x} is given by  
\begin{equation*}
    H_k = \nabla^2_x \operatorname{Tr}(-G \mathcal{A}(x^{(t)}))
\end{equation*}
where \( G \) is the semidefinite matrix solution of Problem~\eqref{eq: SDP-PEP-Generic}. We choose not to include such regularization since it would negatively impact the computational efficiency of our method. The number of variables in \( x \) increases with the number of step sizes to optimize, which, in turn, increases the dimension of the Hessian matrix to consider.  
\end{remark}
\noindent Ideally, we would like to perform the update \( x^{(t+1)} = x^{(t)} + \Delta x^{(t)} \), where \( \Delta x^{(t)} \) is the optimal solution of Problem~\eqref{optPEP_lin}. However, in practice, this step is not always valid, as it does not necessarily ensure a decrease in the objective function. To address this issue, we incorporate a trust-region method into our algorithm. The trust-region approach adaptively bounds the search space of subproblem~\eqref{optPEP_lin} within a safe region—defined by a radius \( D \)—where our linear approximation of the problem remains accurate and the solution of subproblem~\eqref{optPEP_lin} provides an update that decreases the objective function. This is achieved by adding a constraint on the norm of the variable \( \Delta x \):

\begin{equation}\tag{SLM-Subproblem}\label{SSP-sub}
\begin{aligned} 
& \underset{\Delta x}{\text{minimize}}
& & c^{\top} \Delta x  \\ 
& \text{subject to}
& &  \mathcal{A}(x^{(t)}) + \sum \frac{\partial \mathcal{A}}{\partial x_i}(x^{(t)}) \Delta x_i \succeq 0 ,\\
& & & b - d^{\top}(x^{(t)} + \Delta x)= 0,\\
& & &  a_i^{\top} (x^{(t)} + \Delta x) \geq 0 ,\\,
& & & \frac{1}{2} \Delta x^\top \Delta x \leq D.
\end{aligned}
\end{equation}
The next step in the trust-region procedure is to define an update rule for the radius and the current iterate, based on the performance of the update computed by solving Problem~\eqref{SSP-sub}. Let \( D^{(t)} \) and the current iterate \( x^{(t)} = [\tau^{(t)}, \lambda^{(t)}, \alpha^{(t)}] \) be given, along with the update \( \Delta_x^{(t)} = [\Delta \tau^{(t)}, \Delta \lambda^{(t)}, \Delta \alpha^{(t)}] \). To assess the update's performance, we introduce the performance criterion \( P \), defined as the ratio of the actual improvement—represented by the difference in upper bounds on the performances of the first-order algorithms \( \mathcal{M}\left(\alpha_{i,k}^{(t)}\right) \) and \( \mathcal{M}\left(\alpha^{(t)} + \Delta \alpha^{(t)}\right) \)—given by  
\[
\Delta W = w^{\text{dual}}\left(\alpha^{(t)} + \Delta \alpha^{(t)}\right) - w^{\text{dual}}\left(\alpha^{(t)}\right),
\]  
to the expected improvement given by \( \Delta \tau^{(t)}R^2 \): 
\begin{equation*}
    P = \frac{\Delta W}{\Delta \tau^{(t)}R^2}.
\end{equation*}
\begin{remark}
Note that in our case, since the objective of Problem~\eqref{algo-design-generic} is equal to $\tau$ subject to the constraint $\tau \geq 0$, the objective of the linearized subproblem at step $t$ is $\Delta \tau^{(t)}$, under the condition $\tau^{(t)} + \Delta \tau^{(t)} \geq 0$. Since $\tau^{(t)}$ is always positive and $\Delta x = 0$ is always a feasible point of the linearized subproblem, $\Delta \tau^{(t)}$ is always nonpositive, and thus $P \geq 0$ corresponds to a decrease in the upper bound on the worst-case performance ($\Delta W \leq 0$).  
\end{remark}
\noindent If \( P \geq 0.9 \), the update is considered good, and we set \( D^{(t+1)} = 2 D^{(t)} \), \( \alpha_{i,k}^{(t+1)} = \alpha_{i,k}^{(t)} + \Delta \alpha_{i,k}^{(t)} \). If \( P \leq 0.1 \), the update's performance is poor, so we set \( D^{(t+1)} = D^{(t)}/2\) and \( \alpha_{i,k}^{(t+1)} = \alpha_{i,k}^{(t)} \). In all other cases where \( P \) is between $0.1$ and $0.9$, we set \( \alpha_{i,k}^{(t+1)} = \alpha_{i,k}^{(t)} + \Delta \alpha_{i,k}^{(t)} \) while keeping \( D^{(t+1)} = D^{(t)} \). This process is formalized in the following algorithm:  

\begin{algorithm}[H]
    \caption{Trust Region Procedure}\label{alg:trust_reg}
    \begin{algorithmic}[1]
        \State \textbf{Input:} trust-region radius $D^{(t)}$, current iterate $x^{(t)}$, solution $\Delta x^{(t)} = [\Delta\tau^{(t)}, \Delta \lambda^{(t)}_{i,j}, \Delta \alpha^{(t)}_{i,k}]$ of Problem~\eqref{SSP-sub}
        \State \textbf{Output:} $D^{(t+1)}$, $\alpha^{(t+1)}_{i,k}$
        \State Compute $\Delta W = w^{\text{dual}}\left(\alpha^{(t)} + \Delta \alpha^{(t)}\right) - w^{\text{dual}}\left(\alpha^{(t)}\right)$
    
        \State Compute the performance criterion $P = \dfrac{\Delta W}{\Delta\tau^{(t)}R^2}$
        \If{$P \geq 0.9$}
            \State $D^{(t+1)} = 2D^{(t)}$ and $\alpha^{(t+1)}_{i,k} = \alpha^{(t)}_{i,k} + \Delta \alpha^{(t)}_{i,k}$
        \ElsIf{$P \leq 0.1$}
            \State $D^{(t+1)} = \dfrac{D^{(t)}}{2}$ and $\alpha^{(t+1)}_{i,k} = \alpha^{(t)}_{i,k}$
        \Else
            \State $D^{(t+1)} = D^{(t)}$ and $\alpha^{(t+1)}_{i,k} = \alpha^{(t)}_{i,k} + \Delta \alpha^{(t)}_{i,k}$
        \EndIf
    \end{algorithmic}
\end{algorithm}
\noindent We now present our full algorithm as follows:
\begin{algorithm}[H]
\caption{SLM for optimizing first-order algorithms}
\begin{algorithmic}[1]
    \State \textbf{Input:} Initial step-sizes $\alpha^{(0)}_{i,k}$, initial trust-region radius $D^{(0)}$, Number of iterations $T$
    \State \textbf{Output:} $\alpha^{(T)}_{i,k}$
    \For{$t = 1$ to $T$}
    \State Compute $\mathcal{A}(x^{(t)}), \; \lambda^{(t)}_{i,j}, \; \tau^{(t)}$ by solving problem \eqref{dual-PEP-generic}
    \State Form the vector $x^{(t)} = [\tau^{(t)}, \lambda^{(t)}, \alpha^{(t)}]$
    \State Compute $\Delta x^{(t)} = [\Delta \tau^{(t)}, \Delta \lambda^{(t)}, \Delta \alpha^{(t)}] $ by solving~\eqref{SSP-sub}
    \State Compute $D^{(t+1)}$ and $\alpha^{(t+1)}_{i,k}$ using Algorithm~\eqref{alg:trust_reg}.
    \EndFor
\end{algorithmic}
\end{algorithm}

\FloatBarrier
\section{Comparison on the benchmark of Optimizing Memoryless Gradient Descent} \label{sec:bench}
In this section, we compare the efficiency of our three methods with the BNB-PEP approach developed in~\cite{dasgupta2024BNB} for the task of optimizing the step sizes of memoryless gradient descent~\hyperref[alg:MGD]{(MGD)}, defined as follows
\begin{algorithm}[H]
\caption{Memoryless Gradient Descent (MGD)}
\label{alg:MGD}
\begin{algorithmic}[1]
\State \textbf{Input:} function $f$, $x_0 \in \mathbb{R}^d$.
\For{$i = 0$ \textbf{to} $N - 1$}
    \State $x_{i+1} = x_i - \frac{\alpha_i}{L} \nabla f(x_i)$
\EndFor
\end{algorithmic}
\end{algorithm}
\noindent on $L$-smooth convex functions. This fits our framework as~\hyperref[alg:MGD]{(MGD)} is a fixed-step first-order method and $L$-smooth convex functions is an SDP representable class of functions (see Definition~\ref{def:class_func_sdp}) for which a convex PEP can be derived~\cite{taylor2017smooth}. We choose as a performance criterion the final functional accuracy $f(x_N) - f(x_*)$ and as an initial condition ensuring that the worst-case of~\hyperref[alg:MGD]{(MGD)} is bounded $\|x_0 - x_*\| \leq R$.  Thanks to the homogeneity result in~\cite[Section 4.2.5]{taylor2017thesis}, we can set without loss of generality the smoothness constant to \( L = 1 \) and the initial distance to a minimizer to \( R = 1 \) . In the next figure, we provide upper bounds on the worst-case convergence rates for the memoryless gradient descent with the step sizes given in~\cite{dasgupta2024BNB} and with the step sizes obtained using each of our three methods: Alternating Minimization (AM), the First-Order Method (FOM), and our Successive Linearization Method (SLM). Note that the algorithm design approaches discussed here all depend on the total number of iterations $N$ of memoryless gradient descent. The plot shown in the next figure is obtained by running each method for each value of $N$, yielding optimized step sizes specific to that number of iterations. For each case, we compute an upper bound on the performance using PEP, which is then plotted with respect to $N$. We run our methods for a maximum of $T = 1000$ iterations and use a stopping condition where we exit the optimization process when the improvement in the upper bound on the worst-case is less than \( 10^{-7} \), and the norm of the difference between two successive updates of the step sizes is less than \( 10^{-4} \). We initialize our step sizes with the classical constant step-size strategy, setting \( \alpha^{(0)}_i = 1 \). 

\begin{remark}
The BNB-PEP approach consists of two steps. The first step involves solving the algorithm design problem locally by framing it as a nonconvex QCQP problem (similar to our PEP formulations but with substituting the positive semidefinite matrix $G$ using a Cholesky decomposition $G = ZZ^T$ and imposing suitable conditions on the new variable $Z$) and then using the local solution as a warm-start for the Branch-and-Bound procedure. The full BNB-PEP approach including the Branch-and-Bound procedure should converge to global optimum whereas our methods are only local optimization methods. Note that the BNB-PEP approach can also be used as a local optimization heuristic by forgoing the branch and bound procedure (which was done by the authors in~\cite{dasgupta2024BNB} on this benchmark task for a number of iteration higher than $25$).
\end{remark}

\begin{figure}[!htbp]
    \centering
    \includegraphics[width=0.6\linewidth]{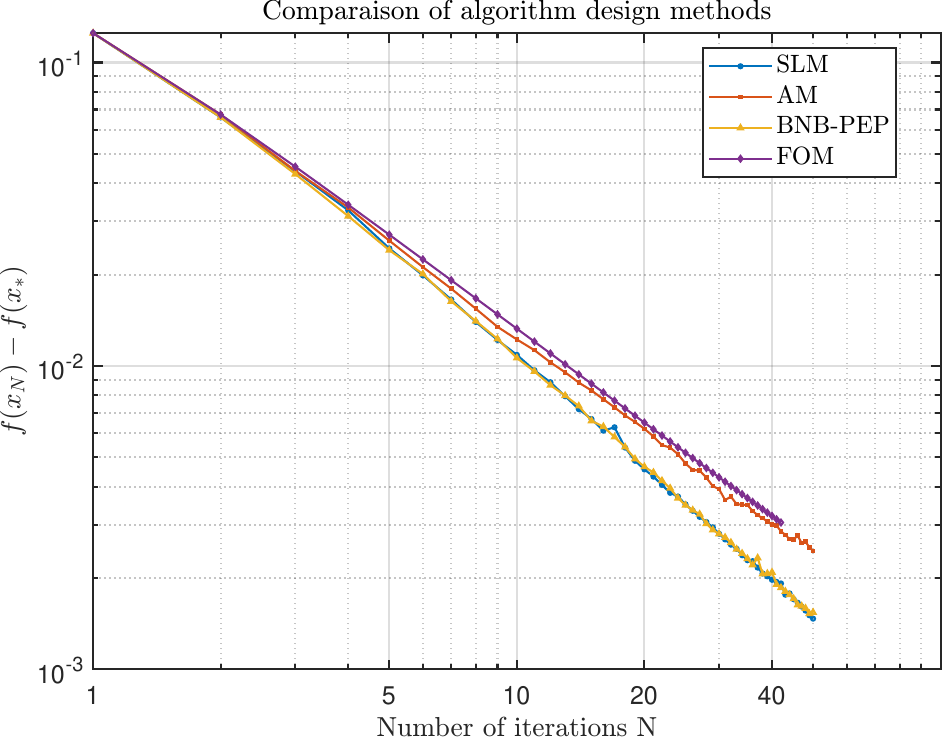}  
    \caption{Worst-case convergence of memoryless gradient descent~\hyperref[alg:MGD]{(MGD)} with the optimized step-sizes given by the BNB-PEP method~\cite{dasgupta2024BNB}, our Alternating Minimization (AM), First-Order method (FOM) and Successive Linearization method (SLM) }
    \label{fig:bench_plot}
\end{figure}
\noindent Figure~\ref{fig:bench_plot} shows that the AM and FOM methods are suboptimal compared to BNB-PEP. As discussed in the previous section, AM can get stuck in nonsmooth regions and fail to converge to an optimum. FOM can also suffer from this issue, as it relies on computing vectors that are gradients of the objective where it is differentiable. Consequently, in nonsmooth regions, the updates do not guarantee a decrease in the objective. Experimentally, we observed that the convergence of these two methods depends significantly on the initialization. For some initializations, the methods may fail to converge to an optimum, while for better initializations (i.e., a good first guess of optimized step sizes) they converge to an optimal point. Moreover, using different stopping criteria that allow these methods to run for more iterations can improve the results. SLM, on the other hand, is competitive with the BNB-PEP approach, even outperforming it in many instances, as shown in Table~\ref{tab:SL-BNB}. Note that the step sizes provided by the BNB-PEP method are claimed to be optimal up to \( N = 25 \). However, Table~\ref{tab:SL-BNB} shows that this is not the case, suggesting that the numerical experiments reported in~\cite{dasgupta2024BNB} using the BNB-PEP approach suffer from numerical issues/bugs due to Gurobi branch-and-bound algorithm's implementation. In Table~\ref{tab:opt_steps_bench}, we present the optimized step-sizes we obtain using SLM for the first few iterations up to $N = 8$.
\begin{table}[!htbp]
    \centering
    \begin{tabular}{ccc@{\qquad}ccc}
    \toprule
     $N$ & SLM & BNB-PEP & $N$ & SLM & BNB-PEP \\
    \midrule
    6  & 0.019895 & 0.020098 & 27 & 0.003181 & 0.003247 \\
    8  & 0.013962 & 0.014055 & 31 & 0.002682 & 0.002722 \\
    9  & 0.012184 & 0.012282 & 32 & 0.002576 & 0.002618 \\
    13 & 0.007946 & 0.007984 & 34 & 0.002379 & 0.002411 \\
    14 & 0.007204 & 0.007376 & 35 & 0.002290 & 0.002318 \\
    16 & 0.006120 & 0.006311 & 37 & 0.002165 & 0.002330 \\
    18 & 0.005382 & 0.005404 & 39 & 0.002027 & 0.002069 \\
    19 & 0.004869 & 0.004941 & 40 & 0.001974 & 0.002089 \\
    20 & 0.004563 & 0.004652 & 43 & 0.001763 & 0.001808 \\
    21 & 0.004318 & 0.004456 & 45 & 0.001703 & 0.001709 \\
    22 & 0.004055 & 0.004185 & 47 & 0.001611 & 0.001616 \\
    23 & 0.003818 & 0.003954 & 48 & 0.001562 & 0.001589 \\
    26 & 0.003333 & 0.003347 & 49 & 0.001508 & 0.001526 \\
    \multicolumn{3}{c}{} & 50 & 0.001470 & 0.001541 \\
    \bottomrule
    \end{tabular}
    \caption{Instances where our Successive Linearization Method (SLM) outperforms BNB-PEP~\cite{dasgupta2024BNB}.}
    \label{tab:SL-BNB}
\end{table}
\begin{table}[!htbp]
\centering
\begin{tabular}{|c|l|}
\hline
\textbf{Number of Iterations $N$} & \textbf{Optimized step-sizes} \\
\hline
1 & [1.500] \\
\hline
2 & [1.414, 1.877] \\
\hline
3 & [1.414, 1.601, 2.189] \\
\hline
4 & [1.414, 1.601, 1.702, 2.459] \\
\hline
5 & [1.414, 1.601, 1.702, 3.526, 1.500] \\
\hline
6 & [1.414, 1.601, 1.702, 4.116, 1.732, 1.500] \\
\hline
7 & [1.414, 1.601, 1.702, 4.552, 1.414, 2.414, 1.500] \\
\hline
8 & [1.414, 1.601, 2.261, 1.414, 5.387, 1.414, 2.414, 1.500] \\
\hline
\end{tabular}
\caption{Optimized step-sizes for memoryless gradient descent obtained using SLM}
\label{tab:opt_steps_bench}
\end{table}
\begin{remark}
    Note that the step sizes in Table~\ref{tab:opt_steps_bench} seem to benefit from a similar structure to one of the step sizes proposed in~\cite{Altschuler2024II,grimmer2024composing} but do not match these schedules. However this similar structure might indicate that these step sizes can be retrieved within the framework proposed in~\cite{grimmer2024composing} of composing smaller step sizes sequences to inductively construct longer ones.
\end{remark}
\noindent We numerically fit the worst-case performance of gradient descent with the SLM step sizes (blue line in Figure~\ref{fig:bench_plot}) using the model \(\frac{1}{\alpha N^{\nu}+\beta}. \) This yields the rate \(\mathcal{O}(N^{-1.237})\), which is slower than the silver-step-size rate \(\mathcal{O}(N^{-1.2716})\)~\cite{Altschuler2024II}, but faster than the BNB-PEP fitted rate \(\mathcal{O}(N^{-1.178})\).  Given the results of Figure~\ref{fig:bench_plot} and Table~\ref{tab:SL-BNB}, we now focus solely on our SLM approach. Next, we compare the computational efficiency of our SLM approach with that of BNB-PEP. In Table~\ref{tab:runtime}, we compare the solver time of the first step of the BNB-PEP procedure (i.e., without the branch and bound procedure, using only the local optimization step on the QCQP) with our SLM approach. We choose IPOPT as the QCQP solver and MOSEK for our SLM approach. We run both methods for a maximum of $T = 1000$ iterations.   
\begin{table}[!htbp] 
    \centering
    \begin{tabular}{|c|c|c|}
    \hline
    Number of iterations $N$ & SLM (s) & BNB-PEP (IPOPT) (s) \\
    \hline
    5  & 0.13  & 10.05  \\
    \hline
    10 & 0.78  & 11.87  \\
    \hline
    15 & 2.06  & 13.05  \\
    \hline
    20 & 7.63  & 61.61  \\
    \hline
    25 & 25.23  & 131.59 (2.19 min) \\
    \hline
    30 & 30.43  & 250.16 (4.17 min) \\
    \hline
    35 & 48.48  & 963.11 (16.05 min) \\
    \hline
    40 & 147.34 (2.46 min)  & 3158.25 (52.64 min)  \\
    \hline
    45 & 174.54 (2.91 min)  & 717.37 (11.96 min)  \\
    \hline
    50 & 371.76 (6.20 min)  & 1294.14 (21.57 min) \\
    \hline
    \end{tabular}
    \caption{Comparison of solver runtime for BNB-PEP using IPOPT for local optimization and our SLM approach.}
    \label{tab:runtime}
\end{table}

\noindent Table~\ref{tab:runtime} indicates that our SLM approach appears to be computationally more efficient than the local optimization (using IPOPT) performed in the BNB-PEP method, especially for larger numbers of step sizes to optimize. Note that both methods in most cases converge before reaching the maximum number of iterations $T = 1000$. In~\cite{dasgupta2024BNB}, the authors provide runtimes for the full BNB-PEP method (including the Branch-and-Bound procedure) on this task. Starting from \( N = 10 \) iterations, the full BNB-PEP method already takes more than one day to run on a supercomputer. To further demonstrate the computational efficiency of our method, we run it respectively for a maximum of optimization steps equal to $T=10$ and $T=100$. We set the maximum iteration parameter of the IPOPT solver to the same values. In the next figure, we provide the worst-case performance rates using the step sizes obtained under these limitations.  

\begin{figure}[!htbp]
    \centering
    \subfloat[$T = 10$ optimization steps.\label{fig:bench_plot_N10}]{
        \includegraphics[width=0.47\textwidth]{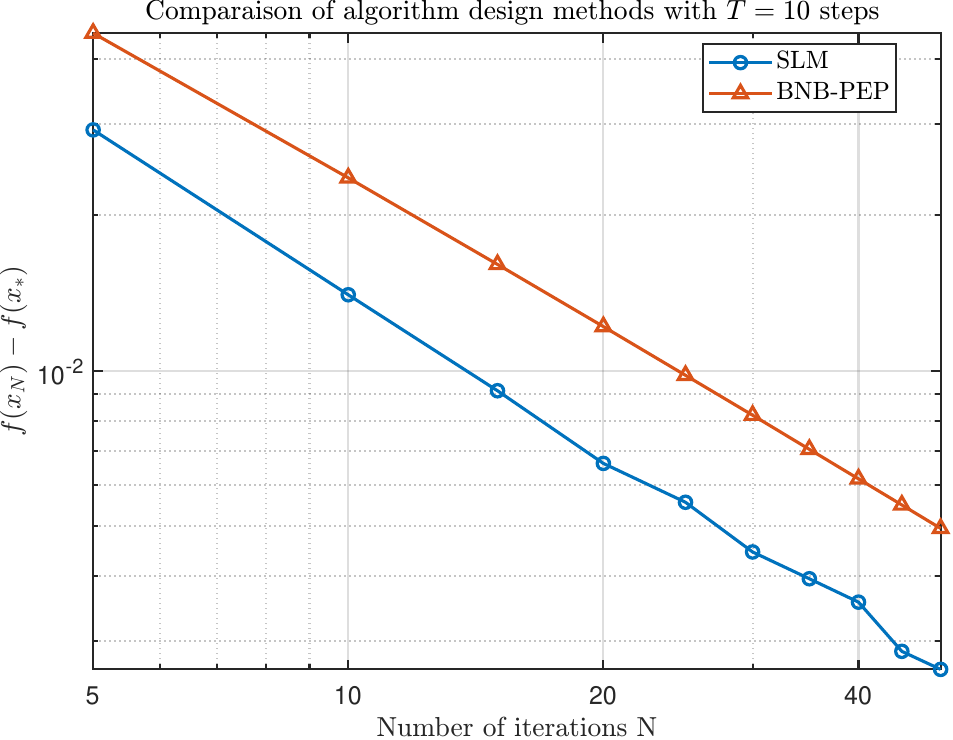}
    }
    \hfill
    \subfloat[$T = 100$ optimization steps.\label{fig:bench_plot_N100}]{
        \includegraphics[width=0.47\textwidth]{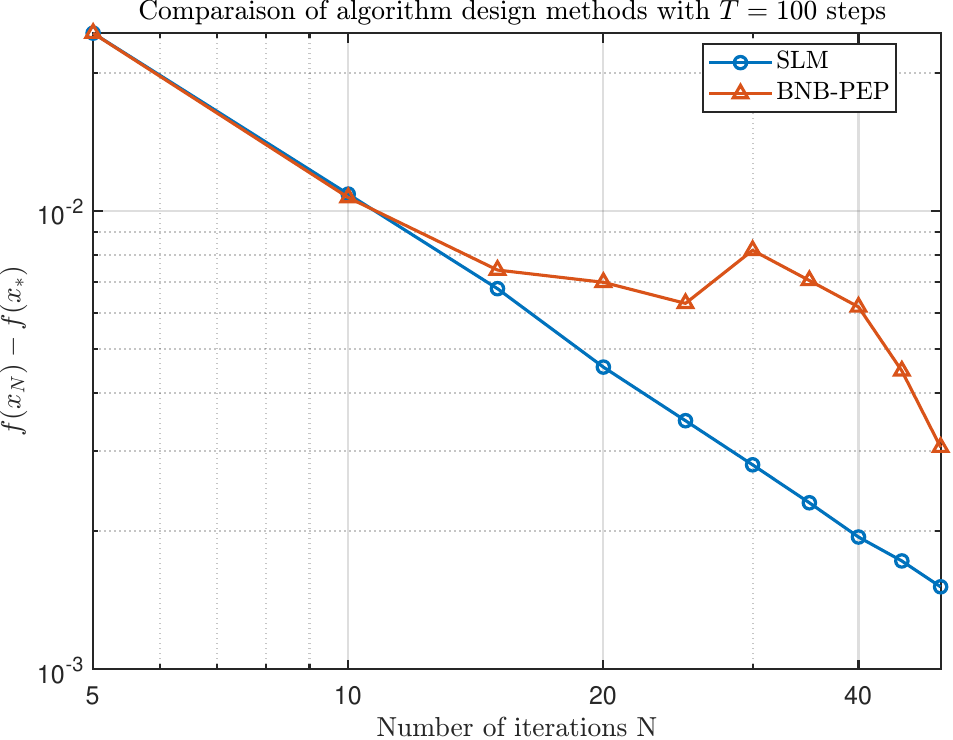}
    }
    \caption{Worst-case convergence of memoryless gradient descent~\hyperref[alg:MGD]{(MGD)} with the optimized step sizes obtained using the BNB-PEP method with the local solver IPOPT, and with the step sizes from our Successive Linearization Method (SLM) under the limits of $T =100$ and $T= 10$ optimization steps.}
    \label{fig:bench_plot_limit}
\end{figure}

\noindent Figure~\ref{fig:bench_plot_limit} shows that, under the constraint of performing few steps to optimize the algorithm ($T = 10$ and $T = 100$ steps), our SLM approach outperforms the BNB-PEP method (with IPOPT). In the next tables, we provide the solver runtimes for performing $T =10$ and $T = 100$ steps for each method. Table~\ref{tab:runtime_comparison_bnb_sl} show that for fewer numbers of iterations ($T = 10$ and $T = 100$) our method has a slightly slower solver run time compared to BNB-PEP. Combined with the numerical results in Figure~\ref{fig:bench_plot_limit} and the solver runtime until convergence in Table~\ref{tab:runtime}, this suggests that our SLM approach has a slower run time per iteration but requires significantly fewer steps to converge.

\begin{table}[!htb]
    \centering
    \small
    \setlength{\tabcolsep}{4pt}
    \renewcommand{\arraystretch}{0.9}
    \begin{tabular}{|c|cc|cc|}
    \hline
    \multirow{2}{*}{$N$}
    & \multicolumn{2}{c|}{$T=10$}
    & \multicolumn{2}{c|}{$T=100$} \\
    \cline{2-5}
    & SLM (s) & BNB-PEP (s)
    & SLM (s) & BNB-PEP (s) \\
    \hline
    5  & 0.11  & 0.024 & 0.17  & 0.10  \\
    10 & 0.17  & 0.10  & 0.95  & 0.64  \\
    15 & 0.54  & 0.31  & 2.66  & 1.14  \\
    20 & 1.26  & 0.86  & 11.11 & 3.07  \\
    25 & 2.43  & 1.94  & 15.91 & 7.95  \\
    30 & 5.43  & 4.32  & 37.81 & 15.60 \\
    35 & 9.14  & 8.90  & 54.63 & 31.22 \\
    40 & 11.47 & 16.73 & 83.68 & 59.12 \\
    45 & 29.60 & 26.93 & 136.97 & 84.67 \\
    50 & 45.76 & 44.68 & 193.35 & 138.66 \\
    \hline
    \end{tabular}
    \caption{Comparison of solver runtimes using SLM and BNB-PEP(IPOPT) for $T=10$ and $T=100$ optimization steps.}
    \label{tab:runtime_comparison_bnb_sl}
\end{table}

\noindent To summarize this section, we presented numerical experiments supporting the competitiveness of our Successive Linearization Method (SLM) with the BNB-PEP approach from~\cite{dasgupta2024BNB} on optimizing step sizes for memoryless gradient descent over smooth convex functions. Moreover, our results suggest that SLM requires fewer steps to converge and performs better under restrictions on the number of steps, indicating overall higher efficiency. In the sequel, we will use SLM to design optimized first-order algorithms in various settings.

\FloatBarrier
\section{Numerical Design of First-Order 
Methods using SLM}\label{sec:opt_algo}
\begin{figure}[!htbp]
    \centering
        \centering
        \includegraphics[width=0.5\textwidth]{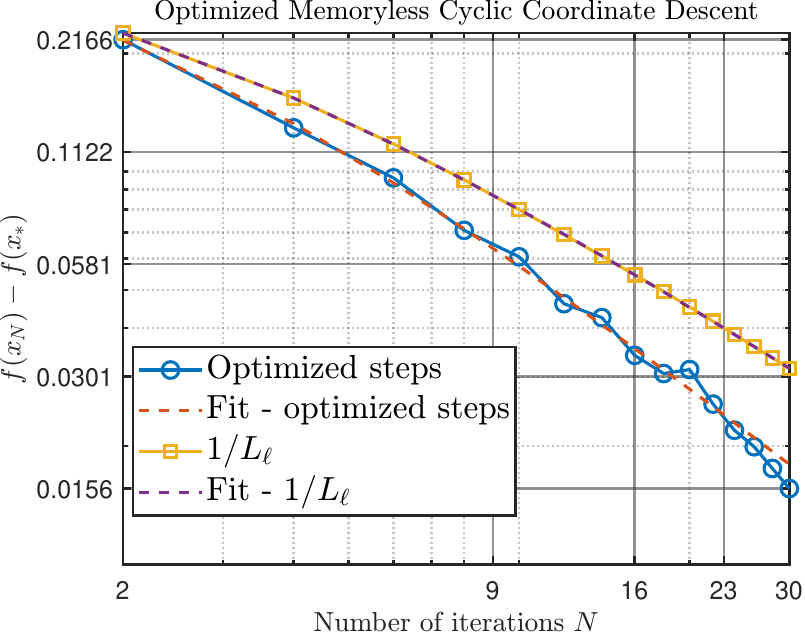}
        \caption{Upper bounds on the worst-case convergence rate of memoryless cyclic coordinate descent with classical steps $\frac{1}{L_\ell}$ (yellow) and with our optimized step-sizes (blue).}
        \label{fig:opt_coords}

\end{figure}

In this section, we apply the SLM method presented above to solve the algorithm design problem and develop fixed-step first-order methods for various convex optimization settings. Numerical results suggest that our approach yields optimized methods that often achieve accelerated convergence rates. For all our experiments, we run our successive linearization method for a maximum of $T=1000$ steps and terminate the optimization process if the improvement in the upper bound on the worst-case performance falls below \(10^{-7} \) and the norm of the difference between successive step-size updates is less than \(10^{-4} \).
\begin{remark}
Note that in this section, the numerical convergence rates presented for the optimized algorithms we design are obtained by solving the Performance Estimation Problem for the algorithm over a range of numbers of iterations $N$, and then fitting a model of the form $\frac{1}{\alpha N^{\nu} + \beta}$ to these data.
\end{remark}
\noindent \textbf{Design of optimized block coordinate-wise algorithms.} We consider the PEP framework derived for block coordinate descent algorithms in~\cite{kamri2025coords} and apply our algorithm design methodology to tune the step sizes of coordinate descent algorithms. The generic framework previously introduced for algorithm optimization extends naturally to this setting. The main difference is that the PEP SDP formulations involve multiple semidefinite matrix variables, which does not hinder the application of the methods presented in this paper. We focus on cyclic algorithms and begin with the subclass of memoryless cyclic coordinate descent algorithms. In order to define block coordinate descent algorithms, we consider a partition of the space \( \mathbb{R}^d \) into \( p \) subspaces:
\begin{equation*}
    \mathbb{R}^d = \mathbb{R}^{d_1} \times \dots \times \mathbb{R}^{d_p},
\end{equation*}
with the corresponding selection matrices \( U_{\ell} \in \mathbb{R}^{d \times d_{\ell}} \), such that:
\begin{equation*}
    (U_1, \dots, U_p) = \mathcal{I}_d,
\end{equation*}
and for any \( x \in \mathbb{R}^d \), we can express \( x \) in terms of its block components as:
\begin{equation*}
    x = (x^{(1)}, \dots, x^{(p)})^T, \quad \text{where} \quad x^{(\ell)} = U_{\ell}^{\top}x \in \mathbb{R}^{d_{\ell}}, \quad \forall \ell \in \{1, \dots, p\},
\end{equation*}
and \( x = \sum_{\ell=1}^{p} U_{\ell} x^{(\ell)} \). Then, memoryless cyclic coordinate descent can be defined as follows.
\begin{algorithm}[H]
\caption{Memoryless Cyclic Coordinate Descent algorithms (MCCD)}\label{alg:MCCD}
\begin{algorithmic}[1]
    \State \textbf{Input:} function $f$ defined over $\mathbb{R}^d$ with $p$ blocks, starting point $x_0 \in \mathbb{R}^d$, number of cycles $K$, and step-sizes $(\alpha_{i})_{i = 1,\dots,pK}$
    \State Define $N = pK$
    \For{$i = 0$ to $N-1$}
        \State Set $\ell = \text{mod}(i,p) + 1$
        \State $x_{i+1} = x_{i} - \frac{\alpha_{i+1} }{L_\ell}U_\ell\nabla^{(\ell)}f(x_{i})$
    \EndFor
\end{algorithmic}
\end{algorithm}
\noindent where \( \nabla^{(\ell)}f(x) = U^{\top}_\ell \nabla f(x) \) is the partial gradient along the \(\ell\)-block of coordinates. We optimize the step sizes \(\{\alpha_{i}\}_{i = 1}^{pK}\) of~\hyperref[alg:MCCD]{(MCCD)}(\(\alpha\)) over the functional class \(\mathcal{F}^{\text{coord}}_{0,\textbf{L}}(\mathbb{R}^d)\), under the assumptions of Setting INIT defined in~\cite{kamri2025coords}, where we only assume that \(\|x_0 - x_*\|_\textbf{L} \leq R\), with \(x_*\) an optimal point of the function. We choose as an initialization for the step sizes $\alpha^{(0)}_i = 1$ for all $i \in \{1,\dots,pK\}$, where $\ell(i) = \text{mod}(i,p) + 1$. By the homogeneity result of~\cite[Theorem 3.1]{kamri2025coords}, we can set $\textbf{L} = (1,\dots,1)$ and the initial distance to a minimizer to $R = 1$ without loss of generality. We further restrict our analysis to two blocks of coordinates ($p = 2$). In Figure~\ref{fig:opt_coords} we compare the numerical upper bounds we obtain using PEP for the classical cyclic coordinate descent algorithm with step sizes $\frac{1}{L_\ell}$ and for the cyclic coordinate descent with our optimized step sizes.
\noindent We also provide numerical fits for these upper bounds. For the classical memoryless coordinate descent algorithm, we obtain a $\mathcal{O}(1/N)$ convergence rate. For the coordinate descent algorithm with optimized step sizes, we have a $\mathcal{O}(1/N^{1.145})$ convergence.

\FloatBarrier \begin{conject}
    Memoryless cyclic coordinate descent with optimized step sizes can achieve a faster convergence rate than $\mathcal{O}(1/N)$. For $2$-blocks algorithms, our experiments suggest a rate of $\mathcal{O}(1/N^{1.145})$.
 \end{conject} \FloatBarrier

\begin{figure}[!htbp]
    \centering
        \centering
        \includegraphics[width=0.5\textwidth]{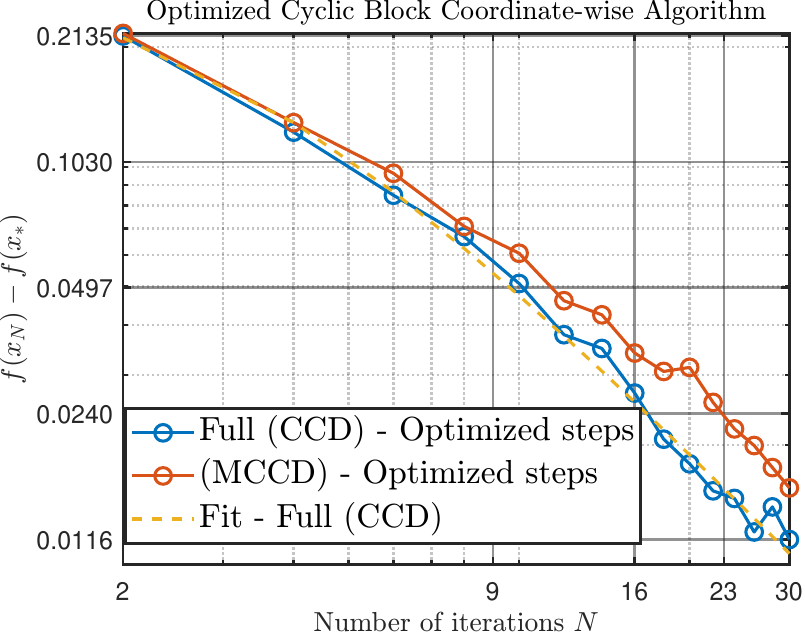}
        \caption{Comparison of upper bounds on the worst-case convergence rate of our optimized full cyclic coordinate descent (blue) and optimized memoryless cyclic coordinate descent (red).}
        \label{fig:opt_coords_mom}
\end{figure}

\noindent We also apply our algorithm design method to the optimization of the step sizes of a full fixed-step cyclic block coordinate-wise algorithm, where we allow the use of all past gradient information to perform the update, similarly to Nesterov's acceleration such that the iterates can be expressed as \begin{equation}\label{eq:coord_algo_fixed}
    x_i = x_0 - \sum_{k = 0}^{i-1} \alpha_{i,k} U_{t(k)}\nabla^{(t(k))} f(x_k),
\end{equation}
with  $t(i) = (i \mod p) + 1$. We place ourselves in the exact same setting as previously for the memoryless coordinate descent. In Figure~\ref{fig:opt_coords_mom}, we compare the upper bounds obtained for both the optimized full cyclic coordinate descent and for the optimized memoryless cyclic coordinate descent highlighting the benefits of using past gradient information. We also provide a numerical fit for the upper bound we have on the convergence of optimized full cyclic coordinate descent with a convergence speed of $\mathcal{O}(1/N^{1.470})$.
\FloatBarrier \begin{conject}
    Full cyclic coordinate descent with optimized step sizes can achieve a faster convergence rate than $\mathcal{O}(1/N)$. For $2$-block algorithms, our experiments suggest a rate of $\mathcal{O}(1/N^{1.47})$.
 \end{conject} \FloatBarrier

\begin{conclusions}
   The numerical upper bounds in Figure~\ref{fig:opt_coords}, along with their corresponding fits, provide strong numerical evidence that accelerated convergence rates for cyclic coordinate descent on coordinate-wise smooth convex functions are achievable without memory (i.e., past gradient information). By accelerated, we mean here $\mathcal{O}(1/N^{\alpha})$ with $\alpha > 1$. This aligns with recent results in~\cite{dasgupta2024BNB,grimmer2024composing,Altschuler2024II} for the gradient descent algorithm on smooth convex functions. The numerical results in Figure~\ref{fig:opt_coords_mom} indicate that using past gradient information, akin to Nesterov's acceleration, helps to further accelerate the convergence of block coordinate-wise algorithms. However, our results suggest an approximate convergence rate of $\mathcal{O}(1/N^{1.47})$, which is significantly slower than the $\mathcal{O}(1/N^{2})$ rate achieved by accelerated random block-coordinate descent algorithms. This supports the comparison made in~\cite{kamri2025coords} regarding the acceleration of random versus cyclic block coordinate descent algorithms. Randomness appears to play a crucial role in accelerating block coordinate-wise algorithms. 
\end{conclusions}

\FloatBarrier
\noindent \textbf{Design of optimized inexact gradient descent algorithms for minimizing the residual gradient norm.} We now consider the optimization of inexact gradient descent algorithms, where updates rely on inexact gradient estimates~\cite{vernimmen2024inex,deklerk2017linesearch}. We consider here relative inexactness on the gradients i.e., we have only access to vectors $d_i$ such that
\begin{equation*}
     \|d_i - \nabla f(x_i)\| \leq \epsilon \|\nabla f(x_i)\|.
\end{equation*}
This arises, for instance, in training machine learning models, where computationally cheaper gradient approximations, such as compressed or truncated gradients, are used. Gradient inexactness can also result from noisy or missing data. We focus on the class of smooth convex functions and begin by analyzing standard memoryless inexact gradient descent algorithms: 
\begin{algorithm}[H]
\caption{Inexact gradient descent (IGD)}\label{alg:IGD}
\begin{algorithmic}[1]
    \State \textbf{Input:} starting Point $x_0$, number of iterations $N$, step-sizes $\alpha_i$, tolerance $\epsilon$
    \For{$i = 0$ to $N-1$}
        \State Compute $d_i$ such that $\|d_i - \nabla f(x_i)\| \leq \epsilon \|\nabla f(x_i)\|$
        \State Update $x_{i+1} = x_i - \frac{\alpha_i}{L} d_i$
    \EndFor
\end{algorithmic}
\end{algorithm}
\noindent Note that inexact gradient descent algorithms can be easily handled within the PEP framework in a similar manner to gradient descent algorithms (see~\cite{taylor2017smooth}). We only need to introduce the extra variables $d_i$ and the conditions $\|d_i - \nabla f(x_i)\| \leq \epsilon \|\nabla f(x_i)\|$. These conditions are readily SDP-representable using the same Gram matrix lifting technique as for gradient descent algorithms, but with the Gram matrix $G$ defined as $G = P^{\top}P$, where $P = [g_0,\dots,g_N,d_0,\dots,d_{N-1},x_0]$. For a more detailed explanation of the PEP formulation for inexact gradient descent, we refer the reader to \cite{deklerk2017linesearch}. Recently, the authors in \cite{vernimmen2024inex} provided the first tight worst-case convergence analysis of inexact gradient descent for smooth convex functions.

\begin{figure}[!htbp] 
    \centering
    \subfloat[$\epsilon = 0.1$\label{fig:opt_inex_grad_01}]{
        \includegraphics[width=0.35\textwidth]{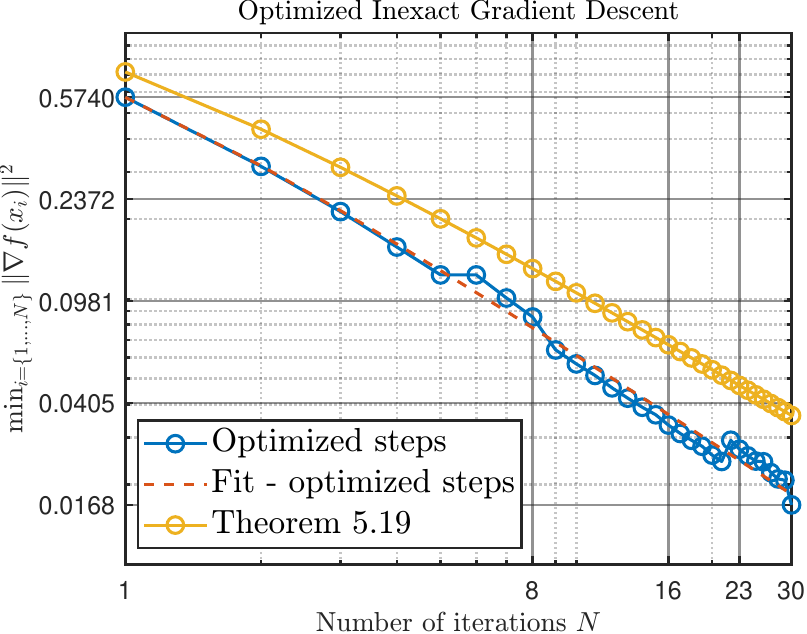}
    }
    \hspace{1em}                         
    \subfloat[$\epsilon = 0.3$\label{fig:opt_inex_grad_03}]{
        \includegraphics[width=0.35\textwidth]{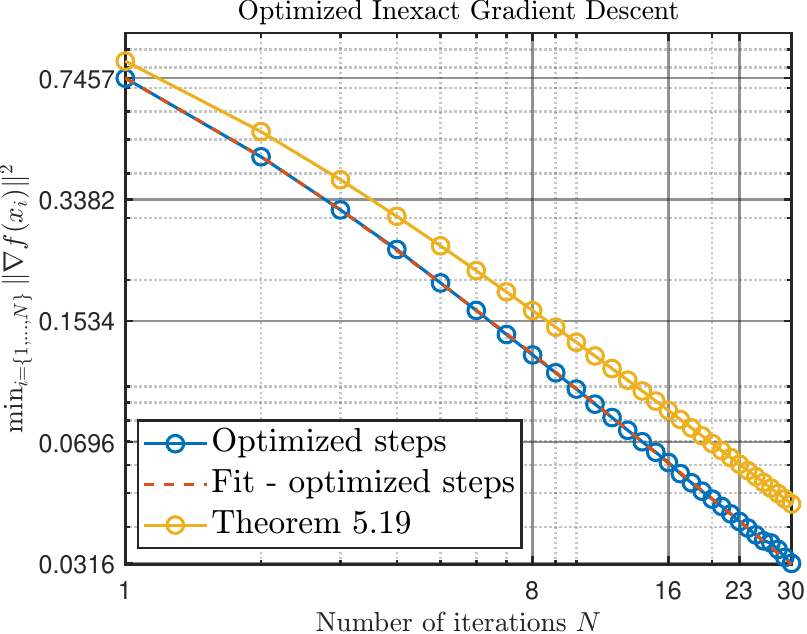}
    }

    \vspace{0.5em}

    \subfloat[$\epsilon = 0.5$\label{fig:opt_inex_grad_05}]{
        \includegraphics[width=0.35\textwidth]{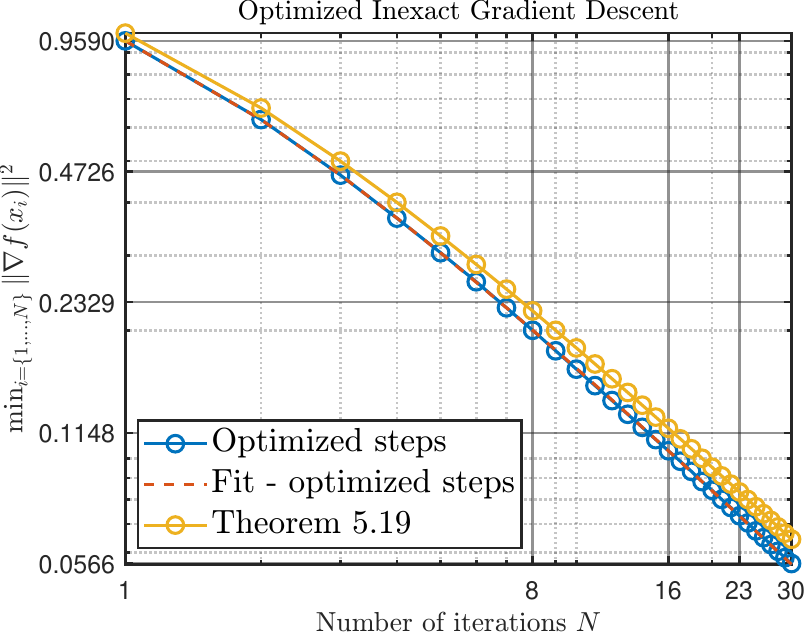}
    }

    \caption{Comparison of worst-case performance of our optimized inexact gradient methods compared to the inexact gradient descent with step-size $\alpha^{(0)}_i = 1/L$ for three different levels of inexactness $\epsilon$.}
    \label{fig:opt_inex_grad}
\end{figure}

\begin{theorem}{\cite[Theorem 2.3]{vernimmen2024inex}}\label{th:inex}
 Algorithm ~\hyperref[alg:IGD]{(IGD)} applied to a convex $L$-smooth function with a constant stepsize $h \in ]0, \frac{3}{2(1 + \epsilon)}]$ and relative inexactness $\epsilon \in ]0, 1[$, started from iterate $x_0$, generates iterates satisfying
\begin{equation*}
    \frac{1}{L} \min_{k \in \{1, \ldots, N\}} \|\nabla f(x_k)\|^2 \leq \frac{f(x_0) - f(x^*)}{\frac{1}{2} + Nh(1 - \epsilon)}.
\end{equation*}
\end{theorem}
\noindent For our numerical experiments, we choose the same performance criterion $\min_{k \in \{1, \ldots, N\}} \|\nabla f(x_k)\|^2$  
and initial condition $f(x_0) - f(x^*) \leq R$ as in Theorem~\ref{th:inex}, and we fix $L = 1$ and $R = 1$. We also apply our algorithm optimization method for several levels of inexactness: $\epsilon = 0.1,\; 0.3,\; 0.5$. We choose $\alpha^{(0)}_i = 1/L$ to initialize the step sizes. In Figure~\ref{fig:opt_inex_grad} we provide numerical values of the worst-case convergence rate of our optimized algorithms, comparing them to the bound of Theorem~\ref{th:inex}. We also performed numerical fits of these upper bounds for our optimized inexact gradient methods which results are summarized in Table~\ref{tab:inex_grad_rates}.

\FloatBarrier \begin{conject}
    Inexact gradient descent with optimized step sizes can achieve a faster convergence rate than $\mathcal{O}(1/N)$ for minimizing the squared gradient norm. We summarize the numerical convergence rates of our optimized inexact gradient descent algorithms in the following table:
\begin{table}[!htbp]
    \centering
    \begin{tabular}{|c|c|}
        \hline
        \textbf{Level of Inexactness $\epsilon$} & \textbf{Worst-Case Convergence Rate} \\ \hline
        0.1                         & $\mathcal{O}(1/N^{1.116})$                          \\ \hline
        0.3                      & $\mathcal{O}(1/N^{1.076})$                           \\ \hline
        0.5                         & $\mathcal{O}(1/N^{1.018})$                     \\ \hline
    \end{tabular}
 \caption{Convergence rates of our optimized inexact gradient descent algorithms with respect to the level of inexactness $\epsilon$.}
\label{tab:inex_grad_rates}
\end{table} 
 \end{conject} \FloatBarrier

\begin{conclusions}
 The numerical worst-case convergence rates of our optimized inexact gradient descent algorithms for levels of inexactness $\epsilon = 0.1, \; 0.3, \; 0.5$ in Figure~\ref{fig:opt_inex_grad}, along with the corresponding numerical fits, indicate that achieving accelerated convergence rates compared to classical inexact gradient descent algorithms with constant step sizes (upper bound of Theorem~\ref{th:inex}) seems possible by carefully tuning the step sizes of the algorithms. Table~\ref{tab:inex_grad_rates}, which summarizes the convergence rates achieved by our optimized algorithms with respect to the inexactness level, suggests that the achievable acceleration depends highly on the level of inexactness. Indeed, it appears that higher acceleration is possible with lower levels of inexactness, whereas for high levels of inexactness, the acceleration tends to become minimal.   
\end{conclusions}

\begin{figure}[!htbp] 
    \centering
    \subfloat[$\epsilon = 0.1$\label{fig:opt_inex_grad_01_mom}]{
        \includegraphics[width=0.35\textwidth]{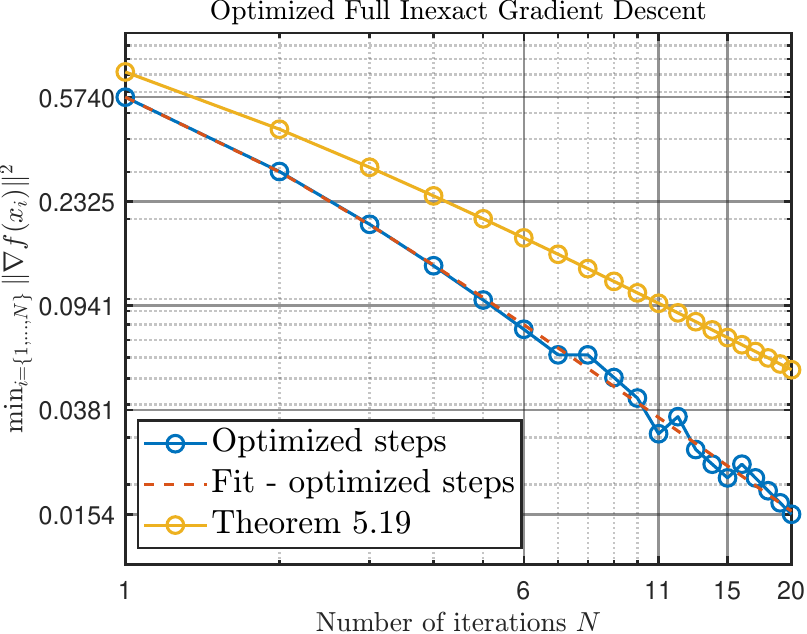}
    }
    \hspace{1em}                         
    \subfloat[$\epsilon = 0.3$\label{fig:opt_inex_grad_03_mom}]{
        \includegraphics[width=0.35\textwidth]{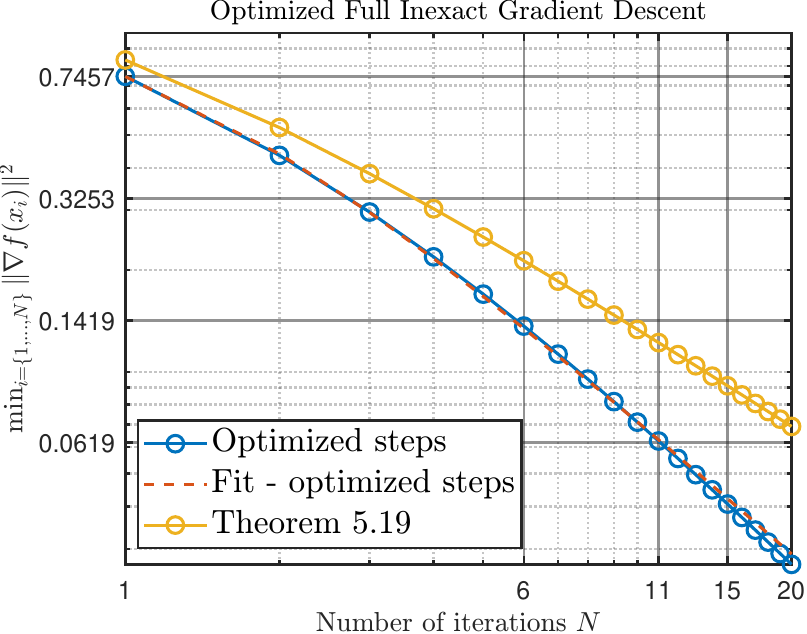}
    }

    \vspace{0.5em}

    \subfloat[$\epsilon = 0.5$\label{fig:opt_inex_grad_05_mom}]{
        \includegraphics[width=0.35\textwidth]{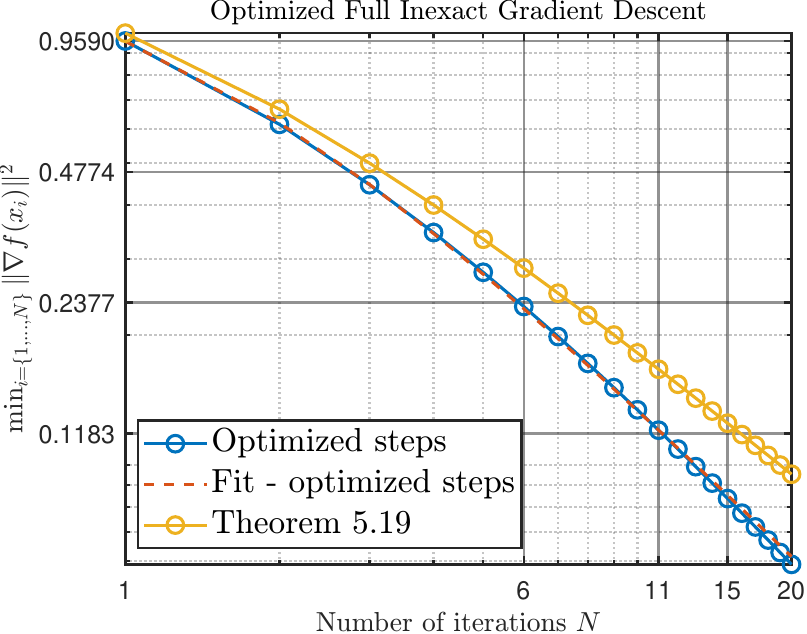}
    }

    \caption{Comparison of upper bounds on the worst-case performance of our optimized full inexact gradient methods compared to the inexact gradient descent with constant step sizes (upper bound from Theorem~\ref{th:inex}) for three different levels of inexactness $\epsilon$.}
    \label{fig:opt_inex_grad_mom}
\end{figure}

\noindent \textbf{Design of full inexact gradient algorithms for minimizing the residual gradient norm.} In an attempt to further accelerate the minimization of the residual gradient norm in the context of inexact optimization of smooth convex functions, we now consider full inexact gradient algorithms, meaning that we allow the use of all past inexact gradient information to perform the updates.
\begin{algorithm}[H]
\caption{Full inexact gradient descent (FIGD)}\label{alg:FIGD}
\begin{algorithmic}[1]
    \State \textbf{Input:} starting Point $x_0$, number of iterations $N$, step-sizes $\alpha_{i,k}$, tolerance $\epsilon$
    \For{$i = 1$ to $N$}
        \State Compute $d_{i-1}$ such that $\|d_{i-1} - \nabla f(x_{i-1})\| \leq \epsilon \|\nabla f(x_{i-1})\|$
        \State Update $x_{i} = x_0 - \sum^{i-1}_{k = 0} \frac{\alpha_{i,k}}{L} d_k$
    \EndFor
\end{algorithmic}
\end{algorithm}
\noindent We place ourselves in the exact same setting as before and make the same assumptions as in the optimization of standard inexact gradient descent algorithms. In Figure~\ref{fig:opt_inex_grad_mom}, we provide numerical upper bounds on the worst-case convergence of our optimized~\hyperref[alg:FIGD]{(FIGD)} and compare them to the bound for the unaccelerated inexact gradient descent (Theorem~\ref{th:inex}). We also performed numerical fits of these upper bounds for our optimized inexact gradient methods whose results are summarized in Table~\ref{tab:inex_grad_rates_mom}.
\FloatBarrier \begin{conject}
Full Inexact gradient descent with optimized step sizes can achieve a faster convergence rate than $\mathcal{O}(1/N)$ for minimizing the squared gradient norm. We summarize the numerical convergence rates of our optimized full inexact gradient descent algorithms in the following table:
\begin{table}[!htbp]
    \centering
    \begin{tabular}{|c|c|c|}
        \hline
        \textbf{Level of Inexactness $\epsilon$} & \textbf{Optimized~\hyperref[alg:FIGD]{(FIGD)}} & \textbf{Optimized~\hyperref[alg:IGD]{(IGD)}} \\ \hline
        0.1 & $\mathcal{O}(1/N^{1.384})$ & $\mathcal{O}(1/N^{1.116})$ \\ \hline
        0.3 & $\mathcal{O}(1/N^{1.342})$ & $\mathcal{O}(1/N^{1.076})$ \\ \hline
        0.5 & $\mathcal{O}(1/N^{1.190})$ & $\mathcal{O}(1/N^{1.018})$ \\ \hline
    \end{tabular}
    \caption{Comparison of the worst-case convergence rates of our optimized full inexact gradient descent and our optimized inexact gradient descent with respect to the level of inexactness $\epsilon$.}
    \label{tab:inex_grad_rates_mom}
\end{table}
 \end{conject} \FloatBarrier

\begin{conclusions}
  The numerical worst-case convergence rates in Figure~\ref{fig:opt_inex_grad_mom}, along with their numerical fits, indicate that the use of memory (past gradient information) indeed clearly helps accelerate the convergence of inexact gradient descent algorithms. However, Table~\ref{tab:inex_grad_rates_mom}, which summarizes the convergence rates achieved by our optimized algorithms, suggests that the achievable rates in the inexact setting are significantly slower than those achieved by optimal gradient descent methods for minimizing the residual gradient norm of smooth convex functions, namely $\mathcal{O}(1/N^2)$ (see~\cite{kim2021OGM-G}).   
\end{conclusions}

\begin{figure}[!htbp] 
    \centering
    \subfloat[$\epsilon = 0.1$\label{fig:opt_inex_01}]{
        \includegraphics[width=0.35\textwidth]{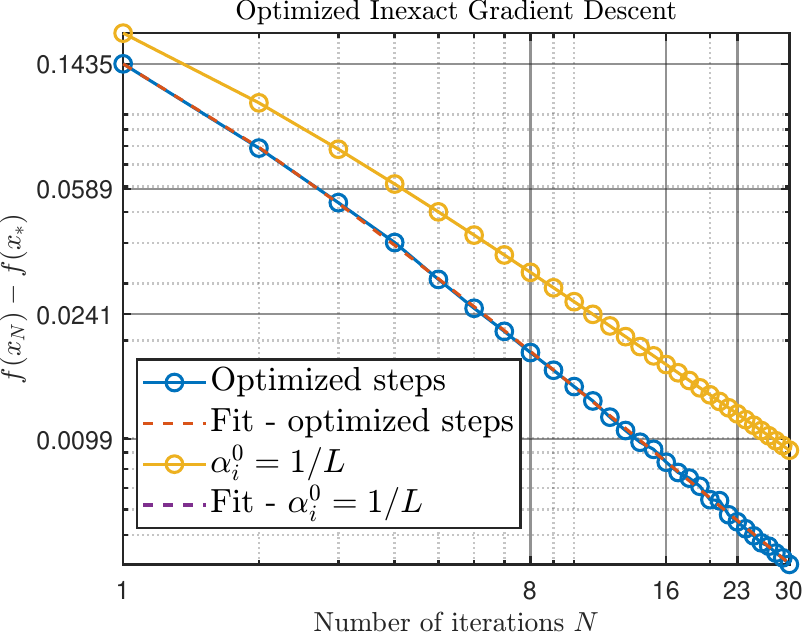}
    }
    \hspace{1em}                         
    \subfloat[$\epsilon = 0.3$\label{fig:opt_inex_03}]{
        \includegraphics[width=0.35\textwidth]{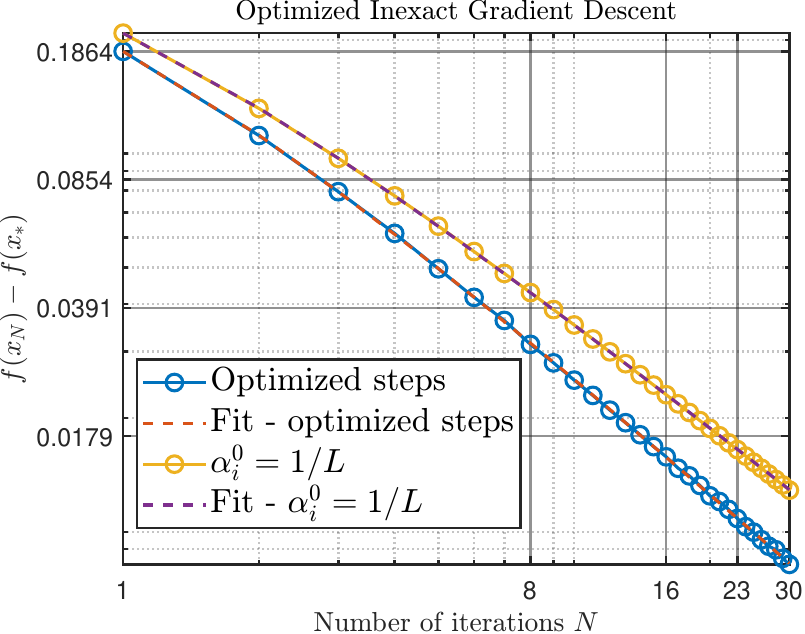}
    }

    \vspace{0.5em}

    \subfloat[$\epsilon = 0.5$\label{fig:opt_inex_05}]{
        \includegraphics[width=0.35\textwidth]{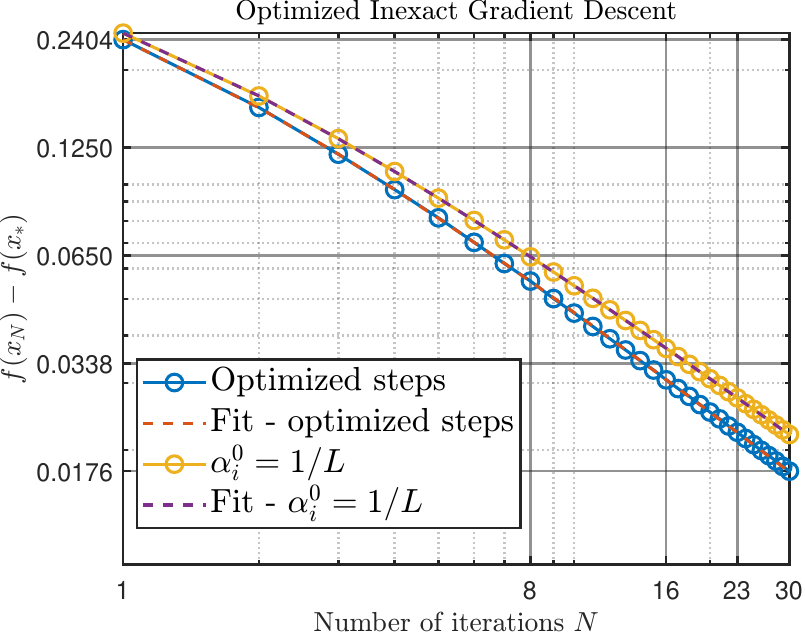}
    }

    \caption{Comparison of the worst-case performance of our optimized inexact gradient methods compared to the inexact gradient descent with constant step sizes $\alpha^{(0)}_i = 1/L$ for three different levels of inexactness $\epsilon$.}
    \label{fig:opt_inex_plot}
\end{figure}

\noindent \textbf{Design of inexact and full (i.e., using past gradient information) inexact gradient descent algorithms for minimizing the functional accuracy.} Previously, we optimized the inexact gradient descent algorithm for minimizing the residual gradient norm on smooth convex functions to compare our results with the only known bound in the inexact setting. However, the residual gradient norm is not the standard performance criterion in convex optimization. Thanks to the flexibility of our PEP-based algorithm design framework, we now optimize inexact gradient descent algorithms for minimizing the functional error $f(x_N) - f(x_*)$, under the initial condition $\|x_0 - x_*\| \leq R$. As before, we fix $L = R = 1$ for our numerical experiments and initialize our optimization with the inexact gradient descent algorithm using a constant step size $\alpha^{(0)}_i = 1/L$. We then compare the performance of our optimized algorithms to this standard algorithm. We begin by optimizing the step sizes of inexact gradient descent without memory and present in Figure~\ref{fig:opt_inex_plot} numerical worst-case convergence rates of our optimized algorithms and the standard inexact gradient descent with step sizes $\alpha^{(0)}_i = \frac{1}{L}$ for levels of inexactness $\epsilon = 0.1,\; 0.3,\; 0.5$. We also performed numerical fits of these upper bounds which results are summarized in  Table~\ref{tab:rates_inex_func}.
 
\begin{figure}[!htbp] 
    \centering
    \subfloat[$\epsilon = 0.1$\label{fig:opt_inex_01_mom}]{
        \includegraphics[width=0.35\textwidth]{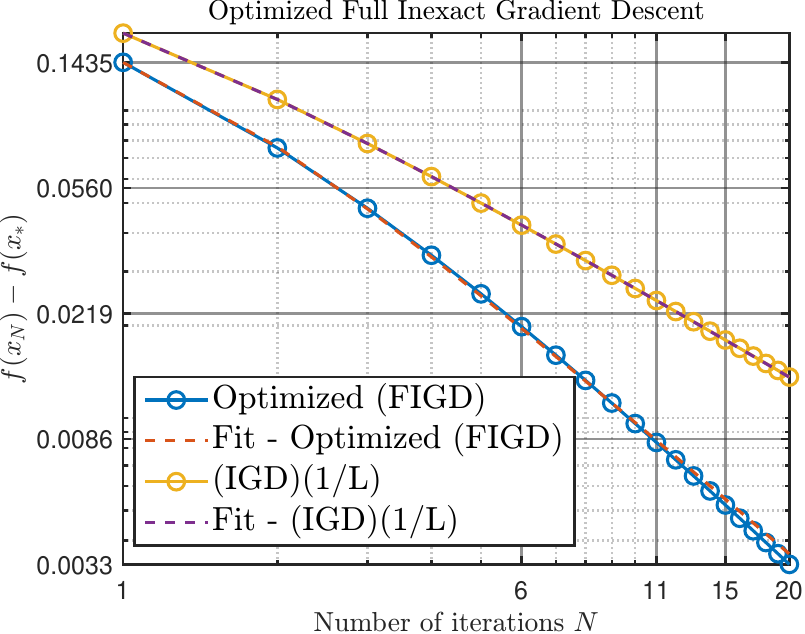}
    }
    \hspace{1em}                         
    \subfloat[$\epsilon = 0.3$\label{fig:opt_inex_03_mom}]{
        \includegraphics[width=0.35\textwidth]{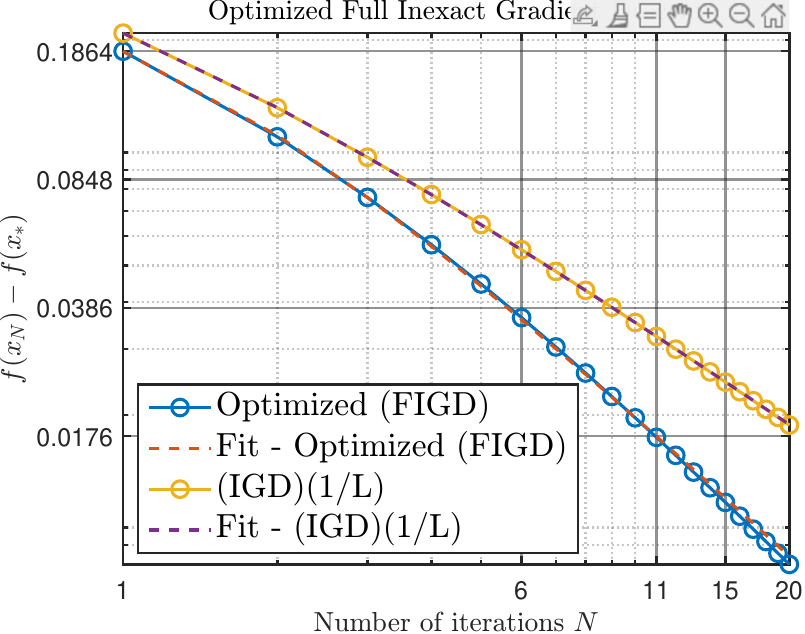}
    }

    \vspace{0.5em}

    \subfloat[$\epsilon = 0.5$\label{fig:opt_inex_05_mom}]{
        \includegraphics[width=0.35\textwidth]{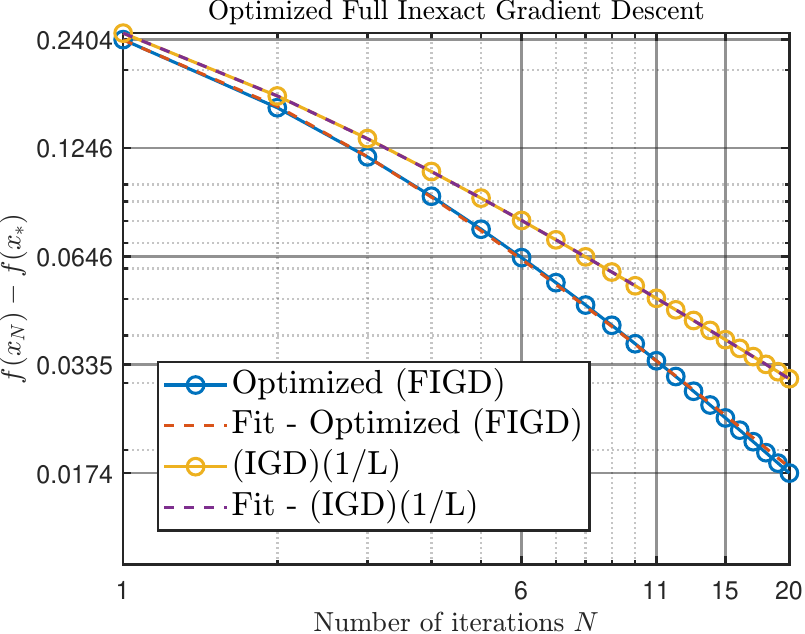}
    }

    \caption{Comparison of worst-case performance of our optimized full inexact gradient methods,~\hyperref[alg:FIGD]{(FIGD)}, compared to the inexact gradient descent with constant step sizes $\alpha^{(0)}_i = 1/L$,~\hyperref[alg:IGD]{(IGD)}, for three different levels of inexactness $\epsilon$.}
    \label{fig:opt_inex_plot_mom}
\end{figure}

\FloatBarrier \begin{conject}
    Inexact gradient descent with optimized step sizes can achieve an accelerated convergence rate compared to inexact gradient descent with step sizes $\frac{1}{L}$ for minimizing the functional accuracy. We summarize the worst-case convergence rates in the following table.  
\begin{table}[!htbp]
    \centering
    \begin{tabular}{|c|c|c|}
        \hline
        \textbf{Level of Inexactness $\epsilon$} & \textbf{\hyperref[alg:IGD]{(IGD)}} & \textbf{Optimized \hyperref[alg:IGD]{(IGD)}} \\
        \hline
        0.1 & $\mathcal{O}(1/N)$ &  $\mathcal{O}(1/N^{1.171})$ \\
        \hline
        0.3 & $\mathcal{O}(1/N^{0.9507})$ & $\mathcal{O}(1/N^{1.052})$  \\
        \hline
        0.5 & $\mathcal{O}(1/N^{0.8786})$ & $\mathcal{O}(1/N^{0.9291})$  \\
        \hline
    \end{tabular}
    \caption{Worst-case convergence rates for our optimized inexact gradient descent algorithms compared to inexact gradient descent with step-sizes $\alpha^{(0)}_i = 1$. }
    \label{tab:rates_inex_func}
\end{table}
 \end{conject} \FloatBarrier
\noindent Note that for $\epsilon = 0.5$, even with optimized steps~\hyperref[alg:IGD]{(IGD)} cannot reach a $\mathcal{O}(1/N)$ rate. We now optimize the full inexact gradient descent algorithm to minimize functional accuracy. Figure~\ref{fig:opt_inex_plot_mom} presents numerical worst-case convergence rates for our optimized algorithms and compares them to those of inexact gradient descent with step sizes $\alpha^{(0)}_i = 1/L$.  We also performed numerical fits of these upper bounds for our optimized algorithm which results are summarized in Table~\ref{tab:rates_inex_func_mom}.
\FloatBarrier \begin{conject}
Full Inexact gradient descent with optimized step sizes can achieve an accelerated convergence rate compared to optimized inexact gradient descent and inexact gradient descent with step sizes $\frac{1}{L}$ for minimizing the functional accuracy. We summarize the worst-case convergence rates in the following table:
\begin{table}[!htbp]
    \centering
    \begin{tabular}{|c|c|c|c|}
        \hline
        \textbf{Inexactness level $\epsilon$} & \textbf{\hyperref[alg:IGD]{(IGD)}} & \textbf{Optimized \hyperref[alg:IGD]{(IGD)}} & \textbf{Optimized \hyperref[alg:FIGD]{(FIGD)}} \\
        \hline
        0.1 & $\mathcal{O}(1/N)$ & $\mathcal{O}(1/N^{1.171})$ & $\mathcal{O}(1/N^{1.450})$ \\
        \hline
        0.3 & $\mathcal{O}(1/N^{0.9507})$ & $\mathcal{O}(1/N^{1.052})$ & $\mathcal{O}(1/N^{1.275})$ \\
        \hline
        0.5 & $\mathcal{O}(1/N^{0.8786})$ & $\mathcal{O}(1/N^{0.9291})$ & $\mathcal{O}(1/N^{1.151})$ \\
        \hline
    \end{tabular}
    \caption{Worst-case convergence rates of the optimized full inexact gradient descent algorithm~\hyperref[alg:FIGD]{(FIGD)} compared to inexact gradient descent with step sizes $\alpha^{(0)}_i = 1/L$,~\hyperref[alg:IGD]{(IGD)} and optimized~\hyperref[alg:IGD]{(IGD)}.}
    \label{tab:rates_inex_func_mom}
\end{table}
 \end{conject} \FloatBarrier

 \begin{conclusions}
    The numerical upper bounds in Figure~\ref{fig:opt_inex_plot} indicate that carefully tuning the step sizes in an inexact gradient descent algorithm accelerates convergence for minimizing functional accuracy. However, convergence remains slower than for minimizing the residual gradient norm. Moreover, the convergence speed, even for standard inexact gradient descent, strongly depends on the level of inexactness in the gradients. The numerical rates in Figure~\ref{fig:opt_inex_plot_mom}, along with the corresponding numerical fits, suggest that incorporating memory further accelerates the convergence of inexact gradient descent for functional accuracy minimization. However, our optimized methods remain significantly slower than optimal gradient methods in the exact minimization setting which have a $\mathcal{O}(1/N^2)$ rate of convergence.
\end{conclusions} 
\FloatBarrier
\noindent \textbf{Design of cyclic gradient descent algorithms in the context of linear convergence.} In the cases where we have linear convergence of gradient descent, our method can be used to design an optimized algorithm compared to the constant step-size gradient descent. Indeed, since we have linear convergence, we can fix a horizon $N_*$, optimize the step sizes for this horizon using our numerical methods, and then consider the cyclic algorithm that repeats the optimized step sizes found by our method. Then if our algorithm verifies after $N_*$ iteration $\|x_{N_*} - x_*\|^2 \leq c(N_*) \|x_0-x_*\|^2$, we have that $\|x_{kN_*} - x_*\|^2 \leq c(N_*)^k \|x_0-x_*\|^2, \; \forall k \in \mathbb{N}$.

\begin{algorithm}[H]
\caption{Cyclic gradient descent}\label{alg:CGD}
\begin{algorithmic}[1]
    \State \textbf{Input:} Starting point $x_0$, number of steps $N$, length of cycles $N_*$, step-sizes $\alpha_1,\dots,\alpha_{N_*}$
    \For{$i = 0$ to $N-1$}
        \State Set $\ell = \text{mod}(i,N_{*}) + 1$
        \State $x_{i+1} = x_{i} - \alpha_{\ell} \nabla f(x_{i})$
    \EndFor
\end{algorithmic}
\end{algorithm}
\begin{remark}
Note that although we restrict our experiments to the class of smooth strongly convex functions, the methodology of designing a cyclic algorithm remains valid in any scenario where linear convergence of first-order algorithms is guaranteed~\cite{Necoara2019,hadi2023lin}.  
\end{remark} 

\begin{figure}[!htbp]
    \centering
    \includegraphics[width=0.6\linewidth]{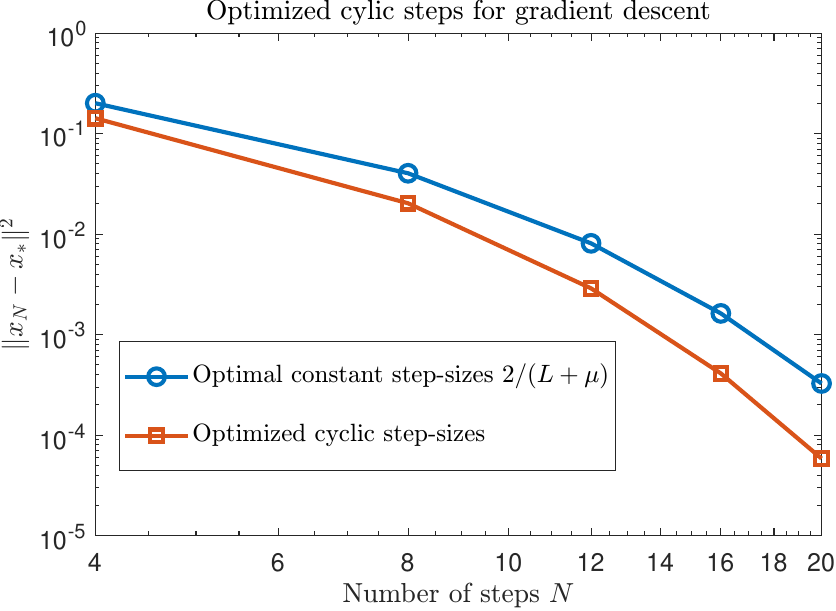}
    \caption{Optimized cyclic steps for gradient descent for $1$-smooth $0.1$-strongly convex functions.}
    \label{fig:cyclic_plot}
\end{figure}

\noindent In Figure~\ref{fig:cyclic_plot}, we fix $N_* = 4$ and compare upper bounds on the worst-case performance of our optimized cyclic gradient descent algorithm to the classical gradient descent algorithm with the optimal constant step size $\frac{2}{L + \mu}$ on $1$-smooth, $0.1$-strongly convex functions.

\begin{conclusions}
  Our optimized cyclic steps improve the worst-case performance of gradient descent on smooth strongly convex functions. To quantify the performance gap, we provide numerical convergence bounds for both methods. For the optimal constant step sizes $\frac{2}{L + \mu}$, we have the worst-case performance $\|x_{4N} - x_*\|^2 \leq 0.2008^N \|x_0 - x_*\|^2, \; \forall N \in \mathbb{N}$. For gradient descent with our cyclic (with cycle length $N_* = 4$) optimized steps, we have the worst-case performance bound$ \|x_{4N} - x_*\|^2 \leq 0.14239^N \|x_0 - x_*\|^2, \; \forall N \in \mathbb{N}$. This provides numerical evidence that our optimized cyclic step-size approach improves the worst-case performance of gradient descent when it converges linearly. Moreover, it is a useful method to trade off between the computational cost of optimizing the steps and the performance of the algorithm. Indeed, larger cycle lengths $N_*$ improve performance but increase the computational cost of optimizing the step sizes.  
\end{conclusions} 

\FloatBarrier
\section{Conclusion}
In this paper, we showed how to formulate the problem of designing fixed-step first-order algorithms in a very generic manner, thanks to the use of the PEP framework. We then provided several methods to optimize first-order algorithms, including a successive linearization method that is competitive with the BNB-PEP method in~\cite{dasgupta2024BNB} on the benchmark task of optimizing memoryless gradient descent on smooth convex functions. Furthermore, we provided numerical evidence indicating the computational efficiency of our method compared to the BNB-PEP approach. We then presented numerically optimized fixed-step first-order algorithms in a range of different settings, showing improved worst-case performance compared to classical fixed-step first-order algorithms available for these optimization settings. Our results suggest that recent findings on the accelerated convergence of first-order algorithms without the use of Nesterov's acceleration over smooth convex functions can be extended to other settings. Open research directions include finding analytical forms for the algorithms presented in this paper, developing more scalable methods for the automatic design of first-order algorithms, as our methods are computationally limited by the cost of solving large-scale SDPs, and extending the design of algorithms using PEP to different classes of algorithms, including, for instance, higher-order methods, decentralized algorithms, and adaptive methods.

\section*{Acknowledgments}
Y. Kamri was supported by the European Union’s MARIE SKŁODOWSKA-CURIE Actions Innovative Training Network (ITN)-ID 861137, TraDE-OPT, and by the FSR program.

\bibliographystyle{plain}  
\bibliography{references}

@article{taylor2017smooth,
  title={Smooth strongly convex interpolation and exact worst-case performance of first-order methods},
  author={Taylor, Adrien B. and Hendrickx, Julien M. and Glineur, Fran{\c{c}}ois},
  journal={Mathematical Programming},
  volume={161},
  number={1-2},
  pages={307--345},
  year={2017},
  publisher={Springer}
}

@phdthesis{taylor2017thesis,
  title     = {Convex Interpolation and Performance Estimation of First-order Methods for Convex Optimization},
  author    = {Taylor, Adrien B.},
  year      = {2017},
  school    = {Université catholique de Louvain},
}

@INPROCEEDINGS{Goujaud2023proofs,
  author={Goujaud, Baptiste and Dieuleveut, Aymeric and Taylor, Adrien},
  booktitle={2023 62nd IEEE Conference on Decision and Control (CDC)}, 
  title={On Fundamental Proof Structures in First-Order Optimization}, 
  year={2023},
  volume={},
  number={},
  pages={3023-3030},
}

@article{drori2014perf,
  author    = {Drori, Yoel and Teboulle, Marc},
  title     = {Performance of First-Order Methods for Smooth Convex Minimization: A Novel Approach},
  journal   = {Mathematical Programming},
  volume    = {145},
  number    = {1–2},
  pages     = {451--482},
  year      = {2014},
  publisher = {Springer}
}

@article{taylor2023ITEM,
  author    = {Adrien Taylor and Yoel Drori},
  title     = {An Optimal Gradient Method for Smooth Strongly Convex Minimization},
  journal   = {Mathematical Programming},
  volume    = {199},
  pages     = {557--594},
  year      = {2023},
}

@inproceedings{toker1995np,
  author    = {Onur Toker and Hitay Ozbay},
  title     = {On the {NP}-Hardness of Solving Bilinear Matrix Inequalities and Simultaneous Stabilization with Static Output Feedback},
  booktitle = {1995 American Control Conference (ACC)},
  pages     = {2525--2526},
  publisher = {IEEE},
  year      = {1995}
}

@article{dasgupta2024BNB,
  author    = {Shuvomoy Das Gupta and Bart P.G. Van Parys and Ernest K. Ryu},
  title     = {Branch-and-Bound Performance Estimation Programming: A Unified Methodology for Constructing Optimal Optimization Methods},
  journal   = {Mathematical Programming},
  volume    = {204},
  pages     = {567--639},
  year      = {2024},
}

@article{kim2016OGM,
  author    = {Donghwan Kim and Jeffrey A. Fessler},
  title     = {Optimized First-Order Methods for Smooth Convex Minimization},
  journal   = {Mathematical Programming},
  volume    = {159},
  pages     = {81--107},
  year      = {2016},
}

@article{nesterov1983FGM,
  author    = {Yurii Nesterov},
  title     = {A Method of Solving a Convex Programming Problem with Convergence Rate $O(1/k^2)$},
  journal   = {Soviet Mathematics Doklady},
  volume    = {27},
  pages     = {372--376},
  year      = {1983}
}

@article{kim2021OGM-G,
  author    = {Donghwan Kim and Jeffrey A. Fessler},
  title     = {Optimizing the Efficiency of First-Order Methods for Decreasing the Gradient of Smooth Convex Functions},
  journal   = {Journal of Optimization Theory and Applications},
  volume    = {188},
  number    = {1},
  pages     = {192--219},
  year      = {2021},
  doi       = {10.1007/s10957-020-01731-7}
}

@article{drori2020GFM,
  author    = {Yoel Drori and Adrien B. Taylor},
  title     = {Efficient First-Order Methods for Convex Minimization: A Constructive Approach},
  journal   = {Mathematical Programming},
  volume    = {184},
  pages     = {183--220},
  year      = {2020},
}

@article{GrimmerLongsteps,
author = {Grimmer, Benjamin},
title = {Provably Faster Gradient Descent via Long Steps},
journal = {SIAM Journal on Optimization},
volume = {34},
number = {3},
pages = {2588-2608},
year = {2024},
}

@misc{Altschuler2023I,
      title={Acceleration by Stepsize Hedging I: Multi-Step Descent and the Silver Stepsize Schedule}, 
      author={Jason M. Altschuler and Pablo A. Parrilo},
      year={2023},
      eprint={2309.07879},
      archivePrefix={arXiv},
}

@Article{Altschuler2024II,
  author    = {Altschuler, Jason M. and Parrilo, Pablo A.},
  title     = {Acceleration by Stepsize Hedging: Silver Stepsize Schedule for Smooth Convex Optimization},
  journal   = {Mathematical Programming},
  year      = {2024},
}

@misc{grimmer2024composing,
      title={Composing Optimized Stepsize Schedules for Gradient Descent}, 
      author={Benjamin Grimmer and Kevin Shu and Alex L. Wang},
      year={2024},
      eprint={2410.16249},
      archivePrefix={arXiv},
      primaryClass={math.OC},
      url={https://arxiv.org/abs/2410.16249}, 
}

@article{nesterov2012Coords,
  author    = {Yurii Nesterov},
  title     = {Efficiency of Coordinate Descent Methods on Huge-Scale Optimization Problems},
  journal   = {SIAM Journal on Optimization},
  volume    = {22},
  number    = {2},
  pages     = {341--362},
  year      = {2012},
  publisher = {Society for Industrial and Applied Mathematics}
}

@article{fercoq2015coord,
  author    = {Olivier Fercoq and Peter Richt{\'a}rik},
  title     = {Accelerated, Parallel and Proximal Coordinate Descent},
  journal   = {SIAM Journal on Optimization},
  volume    = {25},
  number    = {4},
  pages     = {1997--2023},
  year      = {2015},
  publisher = {Society for Industrial and Applied Mathematics},
}

@inproceedings{Boyd1994BMI,
  author    = {S. Boyd and L. El Ghaoui},
  title     = {Alternating Minimization Algorithms for Bilinear Matrix Inequalities},
  booktitle = {Proceedings of the 1994 American Control Conference},
  pages     = {687--691},
  publisher = {IEEE},
  year      = {1994}
}

@book{BonnansShapiro2000,
  author    = {J. Frédéric Bonnans and Alexander Shapiro},
  title     = {Perturbation Analysis of Optimization Problems},
  series    = {Springer Series in Operations Research and Financial Engineering},
  publisher = {Springer},
  year      = {2000},
  isbn      = {978-0-387-98705-7},
}

@article{Correa2004SSP,
  author    = {Rafael Correa and Héctor Ramírez Cabrera},
  title     = {A Global Algorithm for Nonlinear Semidefinite Programming},
  journal   = {SIAM Journal on Optimization},
  volume    = {15},
  number    = {1},
  pages     = {303--318},
  year      = {2004},
}

@article{Fares2002SSP,
  author    = {B. Fares and D. Noll and P. Apkarian},
  title     = {Robust Control via Sequential Semidefinite Programming},
  journal   = {SIAM Journal on Control and Optimization},
  volume    = {40},
  number    = {6},
  pages     = {1791--1820},
  year      = {2002},
}

@book{nemirovsky1983complexity,
  author    = {Nemirovsky, Arkadii and Yudin, David},
  title     = {Problem Complexity and Method Efficiency in Optimization},
  publisher = {Wiley},
  address   = {New York},
  year      = {1983}
}

@article{nemirovsky1984chebyshev,
  author    = {Nemirovsky, A. S. and Polyak, B. T.},
  title     = {The Chebyshev iterative method for solving linear equations},
  journal   = {USSR Computational Mathematics and Mathematical Physics},
  volume    = {24},
  number    = {4},
  pages     = {189--194},
  year      = {1984}
}

@article{polyak1964methods,
  author    = {Polyak, B. T.},
  title     = {Some methods of speeding up the convergence of iteration methods},
  journal   = {USSR Computational Mathematics and Mathematical Physics},
  volume    = {4},
  number    = {5},
  pages     = {1--17},
  year      = {1964}
}

@inproceedings{lin2015universal,
  author    = {Lin, Hongzhou and Mairal, Julien and Harchaoui, Zaid},
  title     = {A Universal Catalyst for First-Order Optimization},
  booktitle = {Advances in Neural Information Processing Systems (NIPS)},
  pages     = {3384--3392},
  year      = {2015}
}

@article{monteiro2013accelerated,
  author    = {Monteiro, Renato D. and Svaiter, Benar F.},
  title     = {An Accelerated Hybrid Proximal Extragradient Method for Convex Optimization},
  journal   = {SIAM Journal on Optimization},
  volume    = {23},
  number    = {2},
  pages     = {1092--1125},
  year      = {2013},
}

@article{beck2009fista,
  author    = {Beck, Amir and Teboulle, Marc},
  title     = {A Fast Iterative Shrinkage-Thresholding Algorithm for Linear Inverse Problems},
  journal   = {SIAM Journal on Imaging Sciences},
  volume    = {2},
  number    = {1},
  pages     = {183--202},
  year      = {2009},
}

@inproceedings{shalev2014accelerated,
  author    = {Shalev-Shwartz, Shai and Zhang, Tong},
  title     = {Accelerated Proximal Stochastic Dual Coordinate Ascent for Regularized Loss Minimization},
  booktitle = {Proceedings of the 31st International Conference on Machine Learning (ICML)},
  pages     = {64--72},
  year      = {2014}
}

@article{allen2017katyusha,
  author    = {Allen-Zhu, Zeyuan},
  title     = {Katyusha: The First Direct Acceleration of Stochastic Gradient Methods},
  journal   = {The Journal of Machine Learning Research (JMLR)},
  volume    = {18},
  number    = {1},
  pages     = {8194--8244},
  year      = {2017}
}

@article{barre2023principled,
  author    = {Barré, Mathieu and Taylor, Adrien B. and Bach, Francis},
  title     = {Principled Analyses and Design of First-Order Methods with Inexact Proximal Operators},
  journal   = {Mathematical Programming},
  volume    = {201},
  number    = {1-2},
  pages     = {185--230},
  year      = {2023},
  doi       = {10.1007/s10107-022-01903-7}
}

@article{jang2023optista,
  author    = {Jang, Uijeong and Das Gupta, Shuvomoy and Ryu, Ernest K.},
  title     = {Computer-Assisted Design of Accelerated Composite Optimization Methods: OptISTA},
  journal   = {arXiv preprint arXiv:2305.15704},
  year      = {2023},
  url       = {https://arxiv.org/abs/2305.15704}
}

@inproceedings{park2022exact,
  author    = {Park, Jongmin and Ryu, Ernest K.},
  title     = {Exact Optimal Accelerated Complexity for Fixed-Point Iterations},
  booktitle = {Proceedings of the 39th International Conference on Machine Learning},
  volume    = {162},
  series    = {Proceedings of Machine Learning Research},
  pages     = {17420--17457},
  publisher = {PMLR},
  year      = {2022},
  month     = {July},
}

@article{goujaud2022optimal,
  author    = {Goujaud, Baptiste and Taylor, Adrien B. and Dieuleveut, Aymeric},
  title     = {Optimal First-Order Methods for Convex Functions with a Quadratic Upper Bound},
  journal   = {arXiv preprint arXiv:2205.15033},
  year      = {2022},
  url       = {https://arxiv.org/abs/2205.15033}
}

@article{hadi2023lin,
  title={Conditions for linear convergence of the gradient method for non-convex optimization},
  author={Abbaszadehpeivasti, H. and de Klerk, E. and Zamani, M.},
  journal={Optimization Letters},
  volume={17},
  pages={1105--1125},
  year={2023},
  publisher={Springer},
  doi={10.1007/s11590-023-01981-2}
}

@article{deklerk2017linesearch,
  title={On the worst-case complexity of the gradient method with exact line search for smooth strongly convex functions},
  author={De Klerk, Etienne and Glineur, François and Taylor, Adrien B},
  journal={Optimization Letters},
  volume={11},
  pages={1185--1199},
  year={2017},
  publisher={Springer}
}

@inproceedings{vernimmen2024inex,
  title={Convergence analysis of an inexact gradient method on smooth convex functions},
  author={Vernimmen, Pierre and Glineur, François},
  booktitle={Proceedings of the European Symposium on Artificial Neural Networks, Computational Intelligence and Machine Learning (ESANN)},
  year={2024},
  address={Bruges, Belgium},
  url={https://www.esann.org/sites/default/files/proceedings/2024/ES2024-171.pdf}
}

@article{Necoara2019,
  author    = {Ion Necoara and Yurii Nesterov and Fran\c{c}ois Glineur},
  title     = {Linear convergence of first order methods for non-strongly convex optimization},
  journal   = {Mathematical Programming},
  volume    = {175},
  pages     = {69--107},
  year      = {2019},
  doi       = {10.1007/s10107-018-1321-1}
}

@mastersthesis{daccache2019performance,
  author       = {Antoine Daccache},
  title        = {Performance Estimation of the Gradient Method with Fixed Arbitrary Step Sizes},
  school       = {Université Catholique de Louvain},
  year         = {2019},
  type         = {Master's thesis}
}

@mastersthesis{eloi2022worst,
  author       = {Diego Eloi},
  title        = {Worst-case Functions for the Gradient Method with Fixed Variable Step Sizes},
  school       = {Université Catholique de Louvain},
  year         = {2022},
  type         = {Master's thesis}
}

@misc{kamri2025coords,
      title={On the Worst-Case Analysis of Cyclic Block Coordinate Descent type Algorithms}, 
      author={Yassine Kamri and François Glineur and Julien M. Hendrickx and Ion Necoara},
      year={2025},
      eprint={2507.16675},
      archivePrefix={arXiv},
      primaryClass={math.OC}, 
}

@phdthesis{kamri2025PhD,
  author = {Kamri, Yassine},
  title = {Contributions to Performance Estimation Problems for the Analysis and Design of First-Order Optimization Methods},
  school = {Universit{\'e} catholique de Louvain},
  year = {2025},
  type = {PhD thesis},
  url = {https://research.dial.uclouvain.be/handle/2078.5/275045}
}

\end{document}